\documentclass[review]{siamart1116}

\usepackage{lipsum}
\usepackage{amsfonts}
\usepackage{graphicx}
\usepackage{epstopdf}
\usepackage{algorithmic}
\ifpdf
  \DeclareGraphicsExtensions{.eps,.pdf,.png,.jpg}
\else
  \DeclareGraphicsExtensions{.eps}
\fi

\numberwithin{theorem}{section}

\newcommand{\TheTitle}{THE STOCHASTIC FITZHUGH-NAGUMO NEURON MODEL IN THE EXCITABLE REGIME EMBEDS A LEAKY INTEGRATE-AND-FIRE MODEL} 
\newcommand{\TheAuthors}{M. E. Yamakou, T. D. Tran, L. H. Duc, and J. Jost}
\newcommand{\TheShortTitle}{FHN neuron model in excitable regime embeds a LIF model} 

\headers{\TheShortTitle}{\TheAuthors}

\title{{\TheTitle}
}

\author{
  Marius E. Yamakou\thanks{Max Planck Institute for Mathematics in the Sciences, 
              Inselstra\ss e 22, D-04103 Leipzig, Germany
    (\email{yamakou@mis.mpg.de},                
              \email{trandat@mis.mpg.de}).}
  \and
  Tat Dat Tran\footnotemark[1]
  \and
  Luu Hoang Duc\footnotemark[1] \thanks{Institute of Mathematics, Viet Nam Academy of Science and Technology, 18 Hoang Quoc Viet Road, 10307 Ha Noi, Viet Nam
               (\email{duc.luu@mis.mpg.de})}
  \and
  J\"urgen Jost\footnotemark[1] \thanks{Santa Fe Institute for the Sciences of Complexity, Santa Fe, NM 87501, USA (\email{jost@mis.mpg.de})
}
}

\usepackage{amsopn}

\usepackage{comment,subfig, float,booktabs, geometry}

\newcommand{\X}{\mathbf X}
\newcommand{\Y}{\mathbf Y}
\newcommand{\Z}{\mathbf Z}

\newcommand{\h}{\mathbf h}

\newcommand{\B}{\mathbf B}

\newcommand{\R}{\mathbb{R}}
\newcommand{\bS}{\mathbf S}
\newcommand{\thetab}{\boldsymbol \theta}

\usepackage{color}

\ifpdf
\hypersetup{
  pdftitle={\TheTitle},
  pdfauthor={\TheAuthors}
}
\fi

\newsiamremark{remark}{Remark}
\numberwithin{equation}{section}

\begin{document}

\maketitle

\begin{abstract}
  In this paper, we provide a complete mathematical construction for a stochastic 
  leaky-integrate-and-fire model (LIF) mimicking the interspike interval (ISI) 
  statistics of a stochastic FitzHugh-Nagumo neuron model (FHN) in the excitable regime, 
  where the unique fixed point is stable. Under specific types of noises, we prove that 
  there exists a global random attractor for the stochastic FHN system. The linearization 
  method is then applied to estimate the firing time and to derive the associated radial 
  equation representing a LIF equation. This result confirms the previous prediction in 
  \cite{Dil12} for the Morris-Lecar neuron model in the bistability regime consisting 
  of a stable fixed point and a stable limit cycle.
\end{abstract}
 


\section{Introduction}\label{sec:intro}
Mathematical modeling has emerged as an important tool to handle the overwhelming structural 
complexity of neuronal processes and to gain a better understanding of their functioning from 
the dynamics of their model equations. However, the mathematical analysis of biophysically 
realistic neuron models such as the 4-dimensional Hodgkin-Huxley (HH) \cite{HH} and the 
2-dimensional Morris-Lecar (ML) \cite{morris81} equations is difficult, as a result of a 
large parameter space, strong nonlinearities, and a high dimensional phase space of the 
model equations. The search for simpler, mathematically tractable (small parameter space, 
weaker nonlinearities, low dimensional phase space) neuron models that still capture all, 
or at least some important dynamical behaviors of biophysical neurons (HH and ML) has been 
an active area of research.

The efforts in this area of research have resulted in easily computable neuron models which 
mimic some of the dynamics of biophysical neuron models. One of the resulting models is  the 
2-dimensional FitzHugh-Nagumo (FHN) neuron model \cite{fitzhugh61}. The FHN model has been 
so successful, because  it is  at the same time mathematically simple and produces a rich 
dynamical behavior that  makes it a model system in many regards, as it reproduces 
the main dynamical features of the HH model. In fact, the HH model has two types of variables, and each type then is combined into a single variable in FHN: 
The ($V,m$) variables of HH correspond to the $v$ variable in FHN, whose fast dynamics represents excitability; 
the ($h,n$) variables correspond to the $w$ variable, whose slow dynamics represents 
accommodation and refractoriness.

The fact that the FHN model is low dimensional makes it possible to visualize the  solution 
and to explain in geometric terms important phenomena related to the excitability and action 
potential generation mechanisms observed in biological neurons. Of course, this comes  
at the expense of  numerical agreement with the biophysical neuron models \cite{yamakou18}. 
The purpose of the model is not a close match with biophysically realistic high dimensional 
models, but rather a mathematical explanation of the essential dynamical mechanism behind the 
firing of a neuron. Moreover, the analysis of such simpler neuron models may lead to the 
discovery of new phenomena, for which we may then  search in the biological neuron models 
and also in experimental preparations.

There is, however, an even simpler model than FHN, the leaky integrate-and-fire model (LIF). 
This is the simplest reasonable neuron model. It only requires a few basic facts about
nerve cells:
they have membranes, they are semipermeable, and they are polarizable.
This suffices to deduce a circuit equivalent to that of  the membrane
potential of the neuron:
a resistor-capacitor circuit. Such circuits charge up slowly when
presented with a current, cross a threshold voltage (a spike), then slowly
discharge.
This behavior is modeled by a  simple 1D equation {together with a reset mechanism}: the leaky
integrate-and-fire neuron model equation \cite{gerstner02}.
Combining sub-threshold dynamics with firing rules has led to a variety
of 1D leaky integrate-and-fire  descriptions of a neuron with a
fixed membrane potential
firing threshold \cite{gerstner02, lansky08}, or with a firing  rate
depending more sensitively on the membrane potential
\cite{pfister06}. In contrast to $n-$dimensional neuron models, $n\geq2$,
such as the HH, ML, and FHN models,
the LIF class of neuron models is less expensive in numerical simulations,
which is an essential advantage when a large network of coupled neurons is considered.

Noise is ubiquitous in neural systems and it may arise from many different sources. 
One  source may come from synaptic noise, that is, the quasi-random release of 
neurotransmitters by synapses or random synaptic input from other neurons. 
As a consequence of synaptic coupling, real neurons operate in the presence of 
synaptic noise. Therefore, most  works in computational neuroscience address 
modifications in neural activity arising from  synaptic noise. Its significance 
can however be judged only if its consequences can be separated from the internal noise, 
generated by the operations of ionic channels \cite{Calvin67}. 
The latter is  channel noise, that is, the random switching of ion channels. 
In many papers  channel noise is assumed to be minimal, because  typically a large number 
of ion channels is involved and fluctuations should  average out, and therefore, the 
effects of synaptic noise should dominate.  Consequently, channel noise is  frequently 
ignored in the mathematical modeling. However, the presence of channel noise can also 
greatly modify the behavior of neurons \cite{white00}. Therefore, in this paper, we 
study the effect of channel noise. Specifically, we add a noise term to the right-hand 
side of the gating equations (the equation for the ionic current variable). 

In the stochastic model, the deterministic fixed point is no longer a solution of 
the system.  The fixed point necessarily needs to vary and adapt to the noise.  To account for this, in the theory of 
random dynamical systems, the notion of a random dynamical attractor was developed as a 
substitute for deterministic attractors in the presence of noise. In the first part of 
this paper, we therefore prove that our system admits a global random attractor, for both additive 
and multiplicative channel noises. This can be seen as a theoretical grounding of our setting.

In \cite{Dil12}, it was shown that a stochastic LIF model constructed with
a radial Ornstein-Uhlenbeck process is embedded in the ML
model (in a bistable regime consisting of a fixed point and limit cycle)
as an integral part of it, closely approximating the sub-threshold
fluctuations of
the ML dynamics. This result suggests that the firing pattern of a
stochastic ML can be recreated using the embedded LIF together with a ML
stochastic firing
mechanism. The LIF model embedded in the ML model captures sub-threshold
dynamics of a combination of the membrane potential and ion channels. 
Therefore, results that can be readily obtained for LIF models can also 
yield insight about ML models. In the second part of this paper, we here address the problem
 to obtain a stochastic LIF model mimicking the interspike interval (ISI) statistics of
 the stochastic FHN model in the excitable regime, where the unique fixed point is stable. 
 Theoretically, we obtain such a LIF model by reducing the 2D FHN model to the one dimensional 
 system that models the distance of the solution to the random attractor as shown in the first 
 part of the paper. In fact, we show that this distance can be approximated to the fixed point, 
 up to a rescaling, as the Euclidean norm $R_t$ of the solution of the linearization of the 
 stochastic FHN equation along the deterministic equilibrium point, and hence the LIF model 
 is approximated by the equation for $R_t$. An action potential (a spike) is produced when 
 $R_t$ exceeds a certain firing threshold $R_t\ge r_0>0$. After firing the process is reset 
 and time is back to zero. The ISI $\tau_0$ is identified with the first-passage time of the 
 threshold, $\tau_0=\inf\{t>0: R_t\geq r_0>0\}$, which then acts as an upper bound of the 
 spiking time $\tau$ of the original system. By defining the firing as a series of first-passage times, 
 the 1D radial process $R_t$ together with a simple firing mechanism based on the detailed FHN model 
 (in the excitable regime), the  firing statistics is shown to reproduce the 2D FHN ISI distribution. 
 We also show that $\tau$ and $\tau_0$ share the same distribution.

The rest of the paper is organized as follows: Sect.~\ref{section2} introduces the deterministic version 
of the FHN neuron model, where we determine the parameter values for which the model is in the excitable 
regime. In Sect.~\ref{section3}, we prove the existence of a global random attractor of the random 
dynamical system generated by the stochastic FHN equation; and furthermore derive a rough estimate for 
the firing time using the linearization method. The corresponding stochastic LIF equation is then 
derived in Sect.~\ref{section4} and its distribution of interspike-intervals is found to numerically match  the stochastic FHN model.

\section{The deterministic model and the excitable regime}\label{section2}

In the
fast time scale $t$, the deterministic FHN neuron model is
\begin{equation}\label{eq:dFHN}
\begin{cases}
dv_t &= (v_t-\displaystyle{\frac{v_t^3}{3}}-w_{t}+I)dt =f(v_t,w_t)dt,\\
dw_t &= \varepsilon (v_t+\alpha-\beta w_t)dt=g(v_t,w_t)dt.
\end{cases}
\end{equation}
where $v_t$ is   the activity of the
membrane potential  and $w_t$ is the recovery current
 that restores the resting state of the model. $I$ is a
constant bias current which can be considered as the effective
external input current.
$0<\varepsilon:=t/\tau \ll1$ is a small singular perturbation
parameter which determines the time scale separation between the
fast $t$ and the slow time scale $\tau$. Thus, the dynamics of $v_t$ is much faster than that of $w_t$. $\alpha$ and $\beta$ are parameters.

\begin{figure}[h!]
\begin{center}
   \includegraphics[width=7cm]{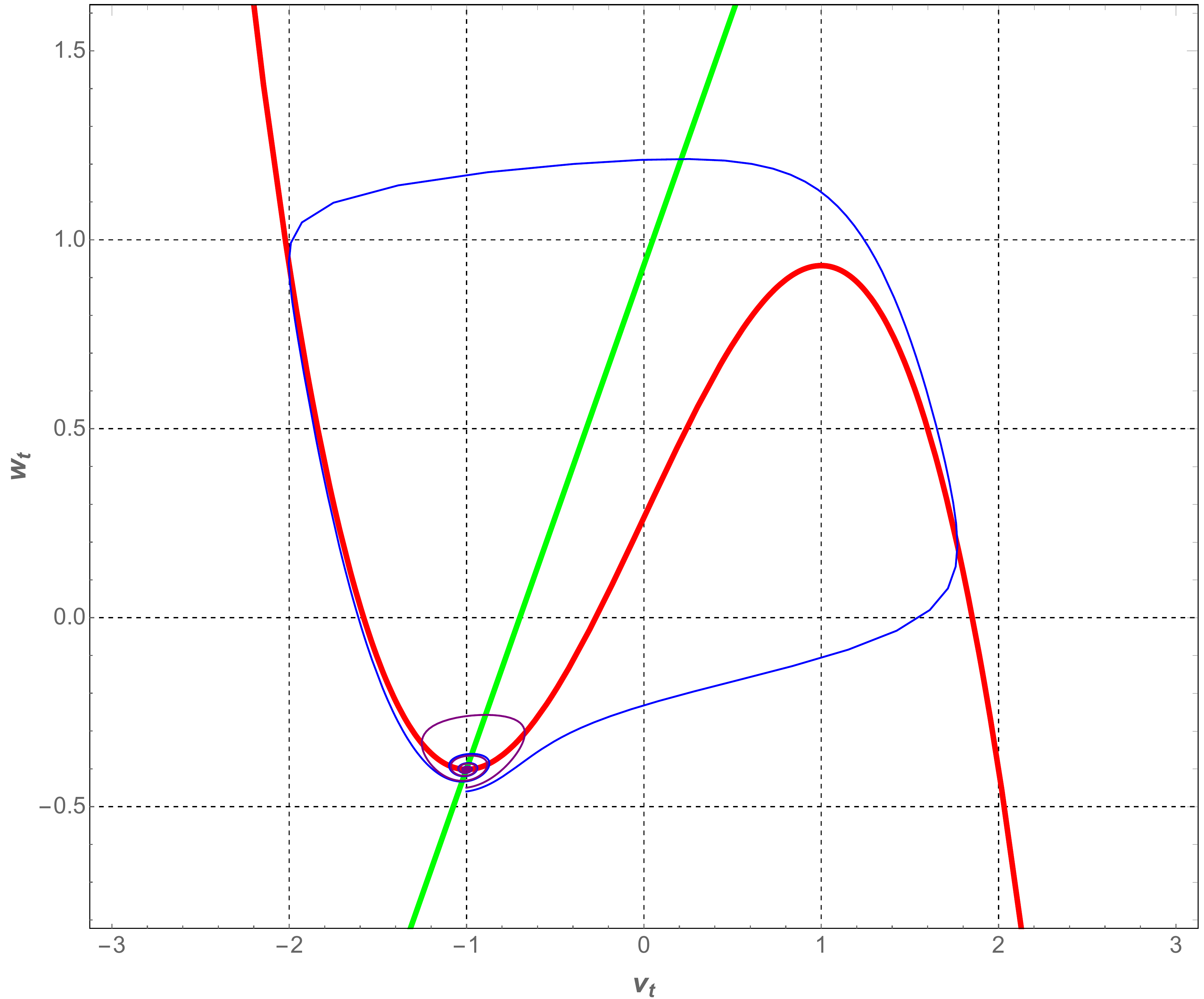}
\end{center}
\caption{\small{
The critical manifold $\mathcal{C}_0$ (red curve) and the $w$-nullcline (green line)
intersect at the unique and stable fixed point $(v_e,w_e)=(-1.00125,-0.401665)$.  Two deterministic trajectories are shown, the purple curve starts at $(-1.00125, -0.45)$ and the blue curve starts at $(-1.00125, -0.46)$. Parameters of system $I = 0.265, \alpha=0.7, \beta=0.75, \varepsilon=0.08$ and the real time for trajectories $T=1000.$}}\label{Fig:1}
\end{figure}

The deterministic critical manifold $\mathcal{C}_0$ defining the
set of equilibria of the \textit{layer problem}
associated to Eq.~\eqref{eq:dFHN} (i.e., the equation obtained from Eq.~\eqref{eq:dFHN} in the singular limit $\epsilon = 0$, see \cite{Kuehn} for a comprehensive introduction to slow-fast analysis), is obtained by solving $f(v,w)=0$ for $w$. Thus, it is given by
\begin{equation}
\mathcal{C}_0=\left\lbrace(v,w)\in\mathbb{R}^2:
w=\displaystyle{v-\frac{v^3}{3}+I}\right\rbrace.
\end{equation}
We note that for Eq.~\eqref{eq:dFHN}, $\mathcal{C}_0$ coincides with the $v$-nullcline (the red curve in Fig.~(\ref{Fig:1})). The stability of points on $\mathcal{C}_0$ as steady
states of the \textit{layer problem} associated to Eq.~\eqref{eq:dFHN} is
determined by the Jacobian scalar $(D_vf)(v,w)=1-v^2$. This shows
that on the critical manifold, points with $|v|>1$ are stable
while points with $|v|<1$ are unstable. It follows that the
branch $v_{-}^*(w)\in (-\infty,-1)$ is stable, $v_0^*(w)\in
(-1,1)$ is unstable, and $v_+^*(w)\in (1,+\infty)$ is stable.

The set of fixed points $(v_e,w_e)$ which define the resting states of the neuron is given by
\begin{equation}
\{(v,w)\in\mathbb{R}^2:f(v,w)=g(v,w)=0\}.
\end{equation}
The sign of the discriminant $\bigtriangleup = (1/\beta-1)^3+\frac{9}{4}(\alpha/\beta-I)^2$,
determines the number of fixed points.
$\mathcal{C}_0$ can therefore intersect the $w$-nullcline ($w=\frac{v+\alpha}{\beta}$) at one,
two or three different fixed points. We assume in this paper that $\bigtriangleup>0$, in which case
we have a unique fixed point given by
\begin{equation}
\begin{cases}
\displaystyle{v_{e}=\sqrt[3]{-\frac{q}{2}-\sqrt{\Delta}}+\sqrt[3]{-\frac{q}{2}+\sqrt{\Delta}}}\\
w_{e}=\frac{1}{\beta}(v_{e}+\alpha).
\end{cases}
\end{equation}

where
\[
p=3\Big(\frac{1}{\beta}-1\Big), \qquad
q=3\Big(\frac{\alpha}{\beta} - I\Big).
\]

Here, we want  to consider the neuron in
the excitable regime \cite{Dil12}. A neuron is  in the excitable regime
when starting  in the basin of
attraction of a unique stable fixed point, {an external pulse} will result into at most
one large excursion (spike) into the phase space
after which the phase trajectory returns back to this fixed point
and stays there \cite{Izhikevich}.

In order to have Eq.~\eqref{eq:dFHN} in the excitable regime, we choose $I, \alpha,$ and $\beta$ such that $\Delta>0$ (i.e., a unique fixed point) and $\varepsilon$ such that the Jacobian (the linearization matrix $M$) of Eq.\eqref{eq:dFHN} at the fixed point $(v_e,w_e)$ has a pair of complex conjugate eigenvalues
\[
-\mu\pm i \nu = \frac{1}{2}(1-v_e^2 - \epsilon \beta) \pm \frac{i}{2} \sqrt{4\epsilon - (1-v_e^2+\epsilon \beta)^2}
\]
with negative real part (i.e., a stable fixed point). In that case, $(v_e,w_e)$ is the only stationary state and there is no limit cycle of system \eqref{eq:dFHN}. In other words, $(v_e,w_e)$ is the global attractor of the system \cite{Izhikevich}. Moreover, to apply the averaging technique \cite{Bax11}, it is necessary that $\mu \ll \nu$, we therefore use through this paper the following parameters of system: $I = 0.265, \alpha=0.7, \beta=0.75, \varepsilon=0.08$ so that $(v_e,w_e) = (-1.00125,-0.401665)$ is the unique stable fixed point and $\frac{\mu}{\nu} = 0.111059 \ll 1$.  Fig.~(\ref{Fig:1}) shows the neuron in the excitable regime. Notice that although every trajectory finally converges to the fixed point, only a small change in the location of the starting point will result in different behavior of the trajectories (see the blue and purple curves).

\section{The stochastic model}\label{section3}

We consider this stochastic FHN model
\begin{equation}\label{eq:1}
\begin{cases}
dv_t &= f(v_t,w_t)dt, \\
dw_t &= g(v_t,w_t)dt+h(w_t)\circ dB_t,
\end{cases}
\end{equation}
where the deterministic fields $f$ and $g$ are given in Eq.~\eqref{eq:dFHN}. There are two important cases:
either $h(w) = \sigma_0$ (additive channel noise) or $h(w)=\sigma_0 w$ (multiplicative channel noise). $\circ dB_t$ stands for the Stratonovich stochastic integral with respect to the Brownian motion $B_t$.\\

\begin{figure}[h!]
\begin{center}
   \includegraphics[width=7cm]{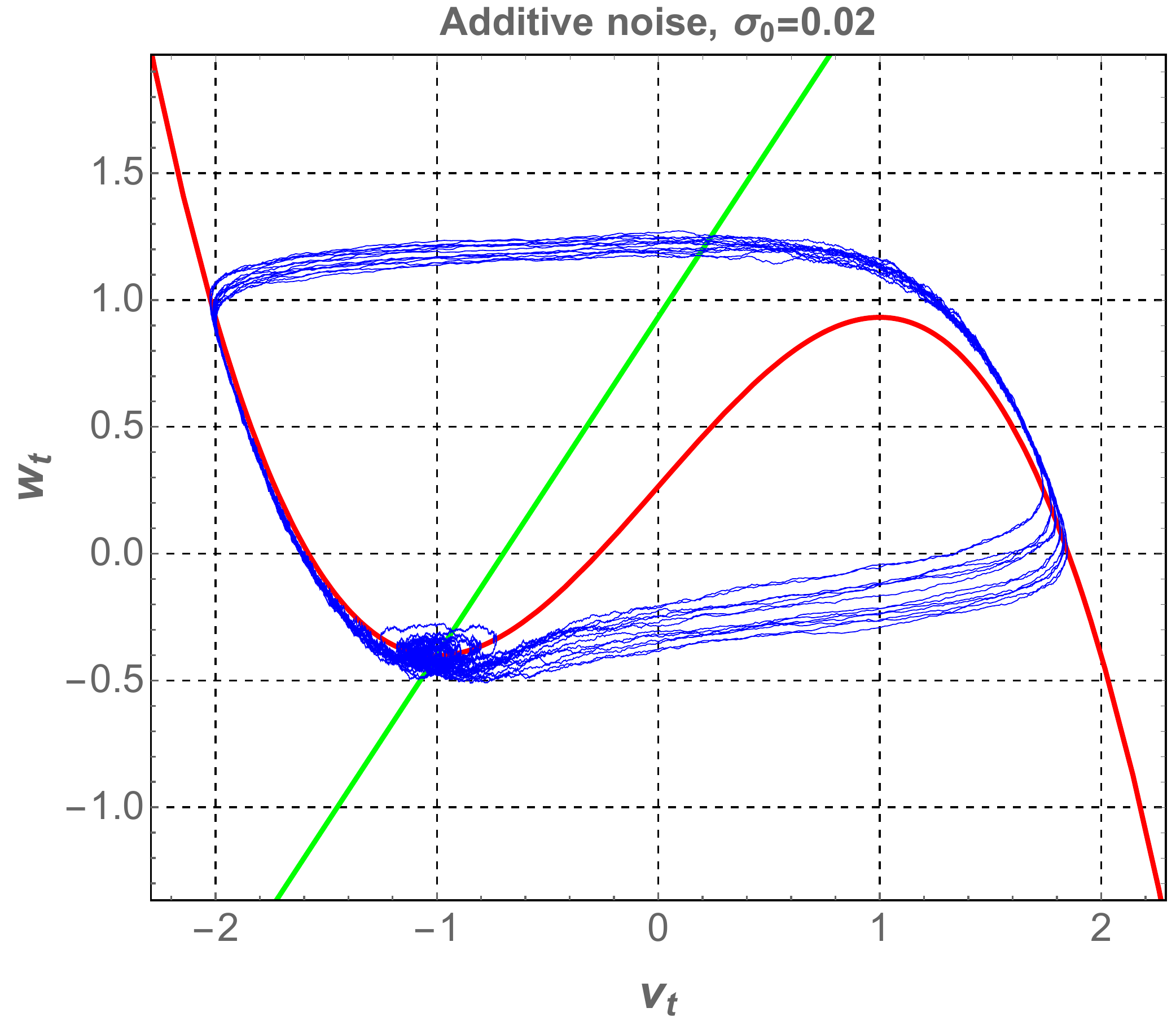} \qquad \includegraphics[width=7cm]{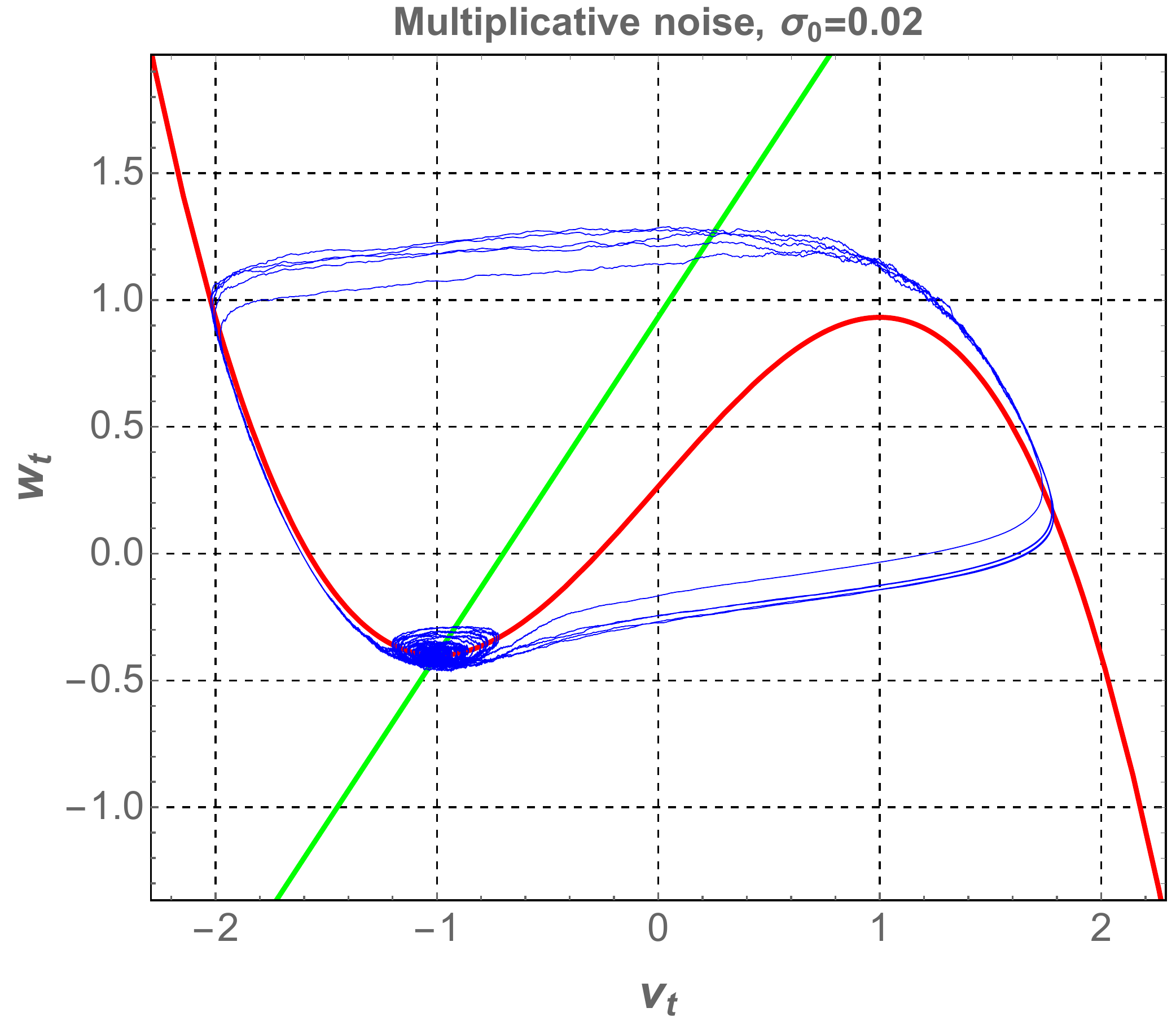}
    \end{center}
\caption{\small Random trajectory of Eq.~\eqref{eq:1} in the excitable regime with chosen parameters of system and the initial condition $(v_0,w_0) = (-1.00125,-0.4)$ for both  additive  and multiplicative noise (we use the StochasticRungeKutta method in {\it Mathematica} with the real time $T=1000$ and the step size $h=0.01$).}
\label{Fig:2}
\end{figure}

\begin{figure}[h!]
\begin{center}
   \includegraphics[width=6cm]{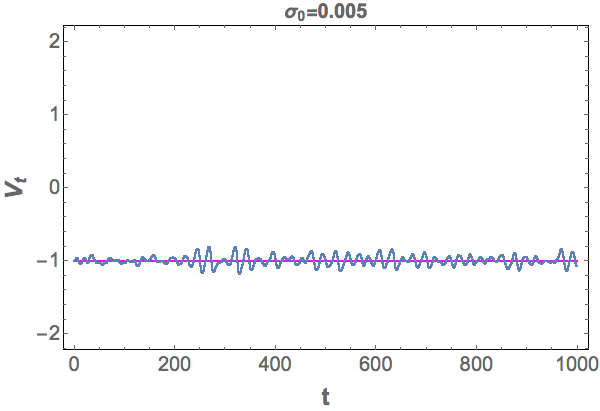}\qquad \includegraphics[width=6cm]{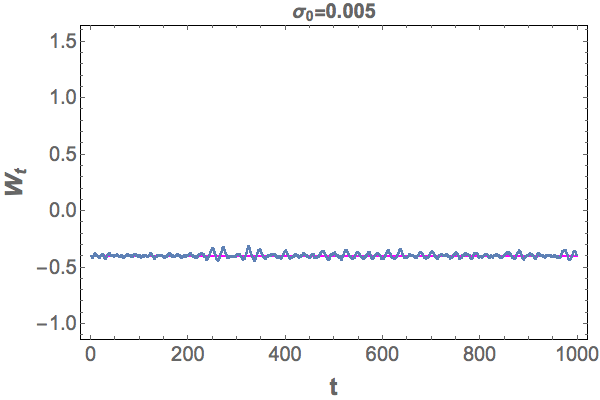}\\
  \includegraphics[width=6cm]{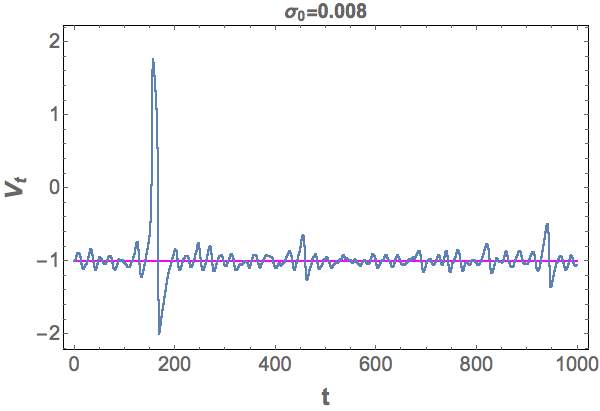}\qquad \includegraphics[width=6cm]{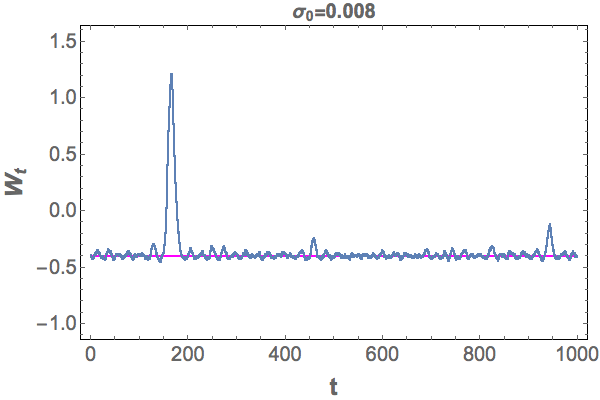}\\
  \includegraphics[width=6cm]{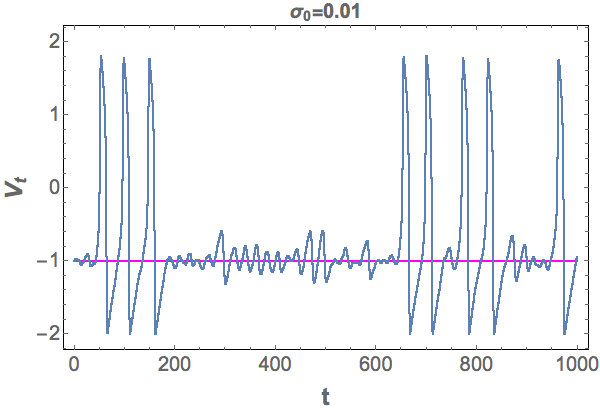}\qquad \includegraphics[width=6cm]{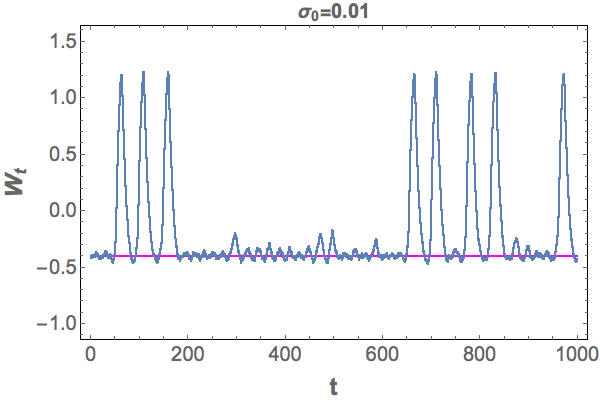}\\
    \includegraphics[width=6cm]{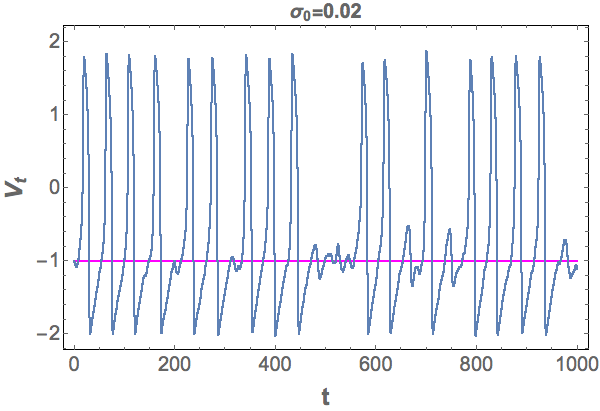}\qquad \includegraphics[width=6cm]{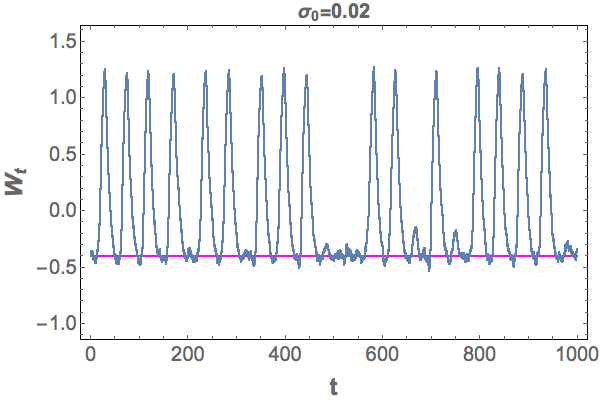}
\end{center}
\caption{\small The components (left column: $V_t$, right column: $W_t$) of a random trajectory of Eq.~\eqref{eq:1} in the excitable regime with chosen parameters of system and the initial condition $(v_0,w_0) = (-1.00125,-0.4)$ for additive noise with $\sigma_0 \in \{0.005, 0.008, 0.01, 0.02\}, T=1000, h=0.01$.}
\label{Fig:3}
\end{figure}

Fig.~\ref{Fig:2} shows the phase portraits of Eq.~\eqref{eq:1} starting with the initial condition $(v_0,w_0) = (-1.00125,-0.4)$, which is in the vicinity of the stable fixed point. Given an initial condition close to the stable fixed point $(v_e,w_e)=(-1.00125,-0.401665)$, the trajectory of the stochastic system might first rotate around the stable fixed point but then {the noise may trigger} a spike, that is, a large excursion into the phase space,  before returning to the neighbourhood of the fixed point; the process repeats itself leading to alternations of small and large oscillations. A similar behavior can be observed when the deterministic system with {an additional} limit cycle is perturbed by noise (as seen in the bistable system \cite{Dil12}).\\
\\
Fig.~\ref{Fig:3} shows that the spiking frequency increases as the amplitude of the noise increases. For a fixed simulation time $T=1000$, the system  spikes  only rarely, if at all,  when the amplitude $\sigma_0 \le 0.005$, but spikes more frequently when $\sigma_0$ increases. This is similar for multiplicative noise.


Let $\X = (v,w)^T$ and $F(\X),H(\X) \in \R^2$ be the drift and  diffusion coefficients of \eqref{eq:1}. The stochastic system is then of the form
\begin{equation}\label{SDE}
d \X_t = F(\X_t)dt + H(\X_t) \circ dB_t,
\end{equation}
where $H(\X) = (0,\sigma_0)^{\rm T}$ for additive noise and 
$
H(\X) = \begin{pmatrix}0&0\\0&\sigma_0 \end{pmatrix} \X =B\X
$
for multiplicative noise. It is easy to check that $F$ is dissipative in the weak sense, i.e.
\begin{eqnarray}\label{dissipativity}
    \langle \X_1 - \X_2,F(\X_1)-F(\X_2) \rangle
    &=& (v_1-v_2)^2 \Big[1- \frac{1}{3} (v_1^2 + v_1 v_2 + v_2^2)\Big]\notag\\
    && -(1-\epsilon)(v_1-v_2)(w_1-w_2) - \epsilon \beta (w_1-w_2)^2\notag\\
    &\leq&  (v_1-v_2)^2 \Big[1 - \frac{1}{12} (v_1-v_2)^2\Big]\notag\\
    &&+ \frac{(1-\epsilon)^2}{2 \epsilon \beta}  |v_1-v_2|^2   + \frac{\epsilon \beta}{2}|w_1-w_2|^2 - \epsilon \beta (w_1-w_2)^2\notag\\
    &\leq& -\frac{1}{12} \Bigg( |v_1-v_2|^2 - 6\Big( 1+ \frac{\epsilon \beta}{2} + \frac{(1-\epsilon)^2}{2 \epsilon \beta}\Big) \Bigg)^2 \\
    & &  +3\Big( 1+ \frac{\epsilon \beta}{2} + \frac{(1-\epsilon)^2}{2 \epsilon \beta}\Big) ^2- \frac{\epsilon \beta}{2} (|v_1-v_2|^2+|w_1-w_2|^2) \notag\\
    &\leq & a- b \|\X_1 - \X_2\|^2
\end{eqnarray}
where
\[
a := 3\Big( 1+ \frac{\epsilon \beta}{2} + \frac{(1-\epsilon)^2}{2 \epsilon \beta}\Big) ^2, \qquad b:= \frac{\epsilon \beta}{2}.
\]
On the other hand, we have
\begin{equation}
|H(\X_1)-H(\X_2)| \le \sigma_0 \Big| w_1-w_2 \Big| \leq \sigma_0 \|\X_1-\X_2\|,
\end{equation}
for multiplicative noise, while $|H(\X_1)-H(\X_2)| \equiv 0$ for additive noise, 
so $H$ is globally Lipschitz continuous.

\subsection{The existence of a random attractor}
In the sequel, we are going to prove that there exists a unique solution $\X(\cdot,\omega,\X_0)$ of \eqref{eq:1} and the solution then generates a so-called {\it random dynamical system} (see e.g. \cite[Chapters 1-2]{arnold}). \\
More precisely, let $(\Omega,\mathcal{F},\mathbb{P})$ be a probability space on which our Brownian motion $B_t$ is defined. In our setting, $\Omega$ can be chosen as $C^0(\R,\R)$, the space of continuous real functions on $\R$ which are zero at zero, equipped with the compact open topology given by the uniform convergence on compact intervals in $\R$, $\mathcal{F}$ as $\mathcal{B}(C^0)$, the associated Borel-$\sigma$-algebra and $\mathbb{P}$ as the Wiener measure. The Brownian motion $B_t$ can then be constructed as the canonical version $B_t(\omega) := \omega(t)$. \\
On this probability space we construct a dynamical system $\theta$ as the Wiener shift
\begin{equation}\label{wienershift}
\theta_t m(\cdot)=m(t+\cdot)-m(t),\quad \forall t\in \R, \forall m \in \bar{\Omega}.
\end{equation}
Then $\theta_t(\cdot): \Omega \to \Omega$ satisfies the group property, i.e.\ $\theta_{t+s} = \theta_t \circ \theta_s$ for all $t,s \in \mathbb{R}$, and is $\mathbb{P}$-preserving, i.e.\ $\mathbb{P}(\theta_t^{-1}(A)) = \mathbb{P}(A)$ for every $A \in \mathcal{F}$, $t \in \R$. The quadruple $((\Omega,\mathcal{F},\mathbb{P},(\theta_t)_{t\in \mathbb{R}})$ is called a {\em metric dynamical system}. \\
Given such a probabilistic setting, Theorem \ref{stadis} below proves that
the solution mapping $\varphi: \R \times \Omega \times \R^2 \to \R^2$ defined by $\varphi(t,\omega)\X_0 := \X(t,\omega,\X_0)$ is a random dynamical system satisfying $\varphi(0,\omega)\X_0 = \X_0$ and the cocycle property
\begin{equation}\label{cocycle}
\varphi(t+s,\omega)\X_0 = \varphi(t,\theta_s \omega)\circ \varphi(s,\omega)\X_0,\qquad \forall t,s \in \R, \omega \in \Omega, \X_0 \in \R^2
\end{equation}

To investigate the asymptotic behavior of the system under the influence of noise, we shall  first check the effect of the noise amplitude on firing.
Under the stochastic scenario, the fixed point $\X_e=(v_e,w_e)$ is no longer the stationary state of the stochastic system \eqref{eq:1}. Instead, we need to find the global asymptotic state as a compact random set $A(\omega) \in \R^2$ depending measurably on $\omega \in \Omega$ such that $A$ is invariant under $\varphi$, i.e. $\varphi(t,\omega)A(\omega) = A(\theta_t \omega)$, and attracts all other compact random sets $D(\omega)$ in the pullback sense, i.e.
\[
\lim \limits_{t\to \infty} d(\varphi(t,\theta_{-t}\omega)D(\theta_{-t}\omega) | A(\omega)) = 0,
\]
where $d(B|A)$ is the Hausdorff semi-distance. Such a structure is called a {\it random attractor} (see e.g. \cite{crauel} or \cite[Chapter 9]{arnold}).\\
The following theorem ensures that the stochastic system~\eqref{eq:1} has a global random pullback attractor. The proof is provided in the Appendix.

\begin{theorem}\label{stadis}
There exists a unique solution of \eqref{SDE} which generates a random dynamical system. Moreover, the system possesses a global random pullback attractor.
\end{theorem}

Theorem~\ref{stadis} shows that every trajectory would in the long run converge to the global random attractor. The structure and the inside dynamics of the global random attractor are still open issues which might help understand the firing mechanism.

\subsection{The normal form at the equilibrium point}
One way to study the dynamics of the stochastic system~\eqref{eq:1} is through its linearization. Therefore, in this section, we shall study the dynamics of \eqref{eq:1} in a small vicinity of the fixed point $\X_e=(v_e,w_e)$. To do that, consider the shift system w.r.t. the fixed point $\X_e$ which has the form
\begin{eqnarray}\label{shifteq}
d(\X_t -\X_e) &=& [F(\X_t) - F(\X_e)] dt + H(\X_t) \circ dB_t \notag\\
&=& \Big[DF(\X_e) (\X_t-\X_e) + \bar{F}(\X_t-\X_e)\Big] dt + H(\X_t) \circ dB_t,
\end{eqnarray}
with initial point $\X_0-\X_e$, where $DF(\X_e)$ is the linearized matrix of $F$ at $\X_e$, $\bar{F}$ is the nonlinear term such that
\[
\begin{aligned}
\|\bar{F}(\X-\X_e)\| &= \Bigg\|\begin{pmatrix} \frac{1}{3}|v+2v_e|(v-v_e)^2 \\ 0 \end{pmatrix}\Bigg\|\\
&\leq \gamma(r) \|\X-\X_e\|,\qquad \forall \|\X-\X_e\| \leq r
\end{aligned}
\]
for an increasing function $\gamma(\cdot): \R_+\to \R_+$, $r\mapsto \frac{r^2}{3} + |v_e|r$, which implies that $\lim \limits_{r \to 0} \gamma(r) =0$. Since $H(\X)$ is either a constant or a linear function, we prove below that system \eqref{shifteq} can be well approximated by its linearized system
\begin{equation}\label{eq:O1}
d{\bar\X}_t =  DF(\X_e) \bar\X_t  dt + H(\bar{\X}_t + \X_e) \circ dB_t,\qquad \bar{\X}_0 = \X_0 -\X_e.
\end{equation}

\begin{theorem}\label{approximation}
    Given $\|\X_0-\X_e\| < r$ and equations \eqref{shifteq}, \eqref{eq:O1}, define the stopping time $\tau = \inf \{t>0: \|\X_t - \X_e\| \geq r\}$. Then there exists a constant $C$ independent of $r$ such that for any $t\geq 0$, the following estimates hold
    \begin{itemize}
        \item For additive noise
        \begin{equation}\label{estthm1}
                \sup_{t\leq \tau}\|\X_t - \X_e -\bar{\X}_t\|\leq C \gamma(r) r.
        \end{equation}

        \item For multiplicative noise
        \begin{equation}\label{estthm}
        E\|\X_{t\wedge\tau} - \X_e -\bar{\X}_{t\wedge\tau}\|^2 \leq C \gamma^2(r) r^2.
        \end{equation}
    \end{itemize}
\end{theorem}
The proof is provided in the Appendix. In practice we  can even approximate \eqref{shifteq} by the following linear system with additive noise
\begin{equation}\label{additivelinear}
d{\tilde\X}_t =  DF(\X_e) \tilde\X_t  dt + H(\X_e) \circ dB_t,\qquad \tilde{\X}_0 = \X_0 -\X_e.
\end{equation}
By the same arguments as in the proof of Theorem~\ref{approximation}, we can prove the following estimate
\begin{equation}\label{estthm2}
    E\|\X_{t\wedge\tau} - \X_e -\tilde{\X}_{t\wedge\tau}\|^2 \leq C r_0^2,
    \end{equation}
for the same stopping time $\tau =  \inf \{t>0: \|\X_t - \X_e\| \geq r_0\}$.

\begin{figure}[h!]
\begin{center}
    \includegraphics[width=6cm]{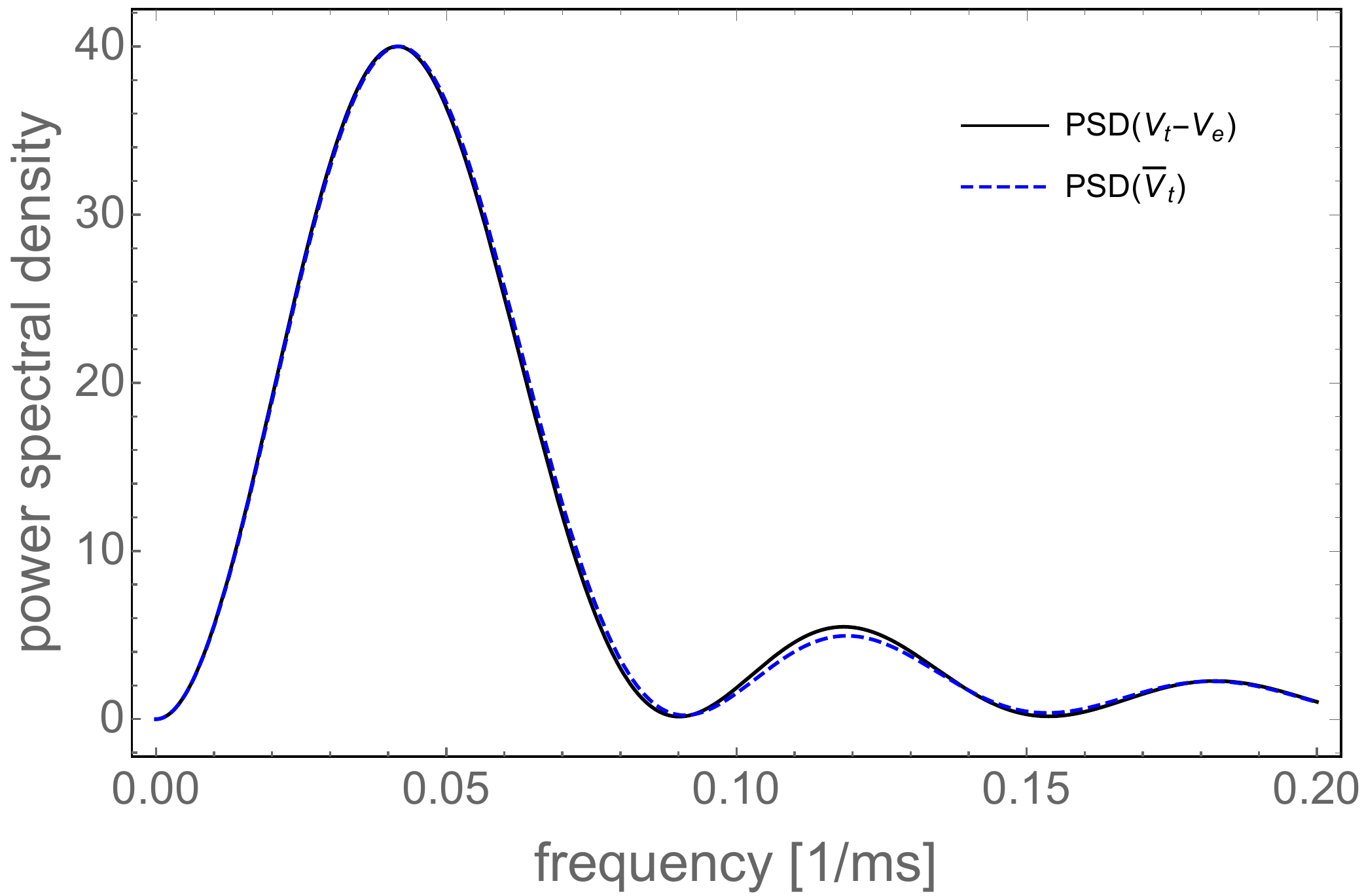}\qquad \includegraphics[width=6cm]{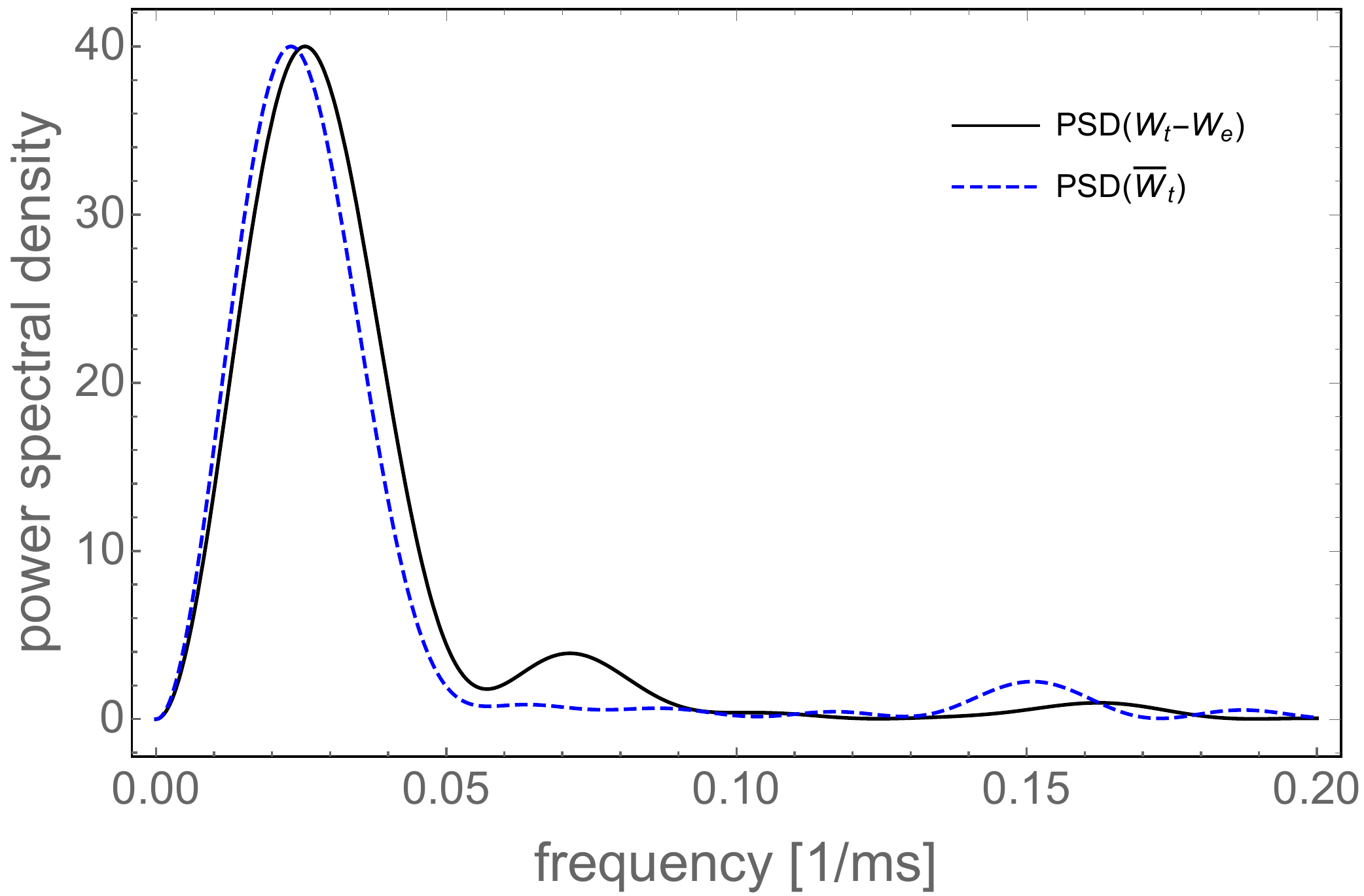}
\end{center}
\caption{\small The process $\{\X_t -\X_e\}_t$ \eqref{shifteq} and the process $\{\bar \X_t\}_t$~\eqref{eq:O1} with the chosen parameters of system, $\sigma_0=0.01$ and the same starting point $(v_0-v_e,w_0-w_e)$ are compared by using the power spectral density. Their spectrum densities are well approximated.}
\label{Fig:5}
\end{figure}

Another comparison between the processes $\{\X_t -\X_e\}_t$ and $\{\bar \X_t\}_t$ can be obtained by using power spectral density estimation (see, for example, \cite[Chapter 7]{Fan2003}). In Fig.~\ref{Fig:5}, the estimated spectral densities of the shifted original and {the linearized process} are plotted. The spectral densities are estimated from paths started from $0$ to $50$ ms of subthreshold fluctuations, and scaled to have the same maximum at $40$.

\section{The embedded LIF model}\label{section4}

In this section, we present two constructive methods to obtain
$1$-D LIF models corresponding to the stochastic FHN in the
excitable regime in Eq.~\eqref{eq:1}. The first method follows \cite{Bax11} (see also
\cite{Dil12}) by constructing the so-called {\it radial
Ornstein-Uhlenbeck equation}. More precisely, we rewrite the
linearized system~\eqref{eq:O1} in the form
\begin{equation}\label{additive}
d\bar{\X}_t = M \bar{\X}_t dt + \begin{pmatrix}0 & 0\\ 0 & \sigma_0\end{pmatrix} d\B_t,
\end{equation}
where $M= DF(\X_e)$ and $\B_t = \begin{pmatrix}B'_t\\ B_t\end{pmatrix}$  is a $2$-D standard Brownian motion. For chosen parameters, $M$ has a pair of complex conjugate eigenvalues $-\mu\pm i \nu$ with $\mu = 0.0312496, \nu = 0.281378$. By transformation $\bar{\Y}_t = Q^{-1} \bar \X_t$ with $Q = \begin{pmatrix}-\nu & m_{11}+\mu\\ 0 & m_{21}\end{pmatrix}$ we obtain
\begin{equation}\label{additive1}
d \bar{\Y}_t = A\bar{\Y}_t dt + C d\B_t,
\end{equation}
where
\[
A = \begin{pmatrix}-\mu & \nu\\ -\nu & -\mu\end{pmatrix} =  \begin{pmatrix}-0.0312496 & 0.281378\\ -0.281378 & -0.0312496\end{pmatrix} ;
\]
\[
C = Q^{-1} \begin{pmatrix}0 & 0\\ 0 & \sigma_0\end{pmatrix}.
\]

We note that $\frac{\mu}{\nu} = 0.111059 \ll 1$, therefore, by applying the technique of time average from \cite[Theorem 1]{Bax11},
$\bar \Y_t$ can be approximated by an Ornstein-Uhlenbeck process up to a rotation, i.e.
\[
\bar \Y_t \sim \bar \Y^{app}_t := \frac{\sigma}{\sqrt{\mu}} Rot_{-\nu t} \bar \bS_{\mu t},
\]
where
$\sigma = \sqrt{\frac{1}{2} \mathrm{tr} (C C^*)} = \sqrt{\frac{-m_{12}}{2 \nu^2 m_{21}}}\sigma_0$, the rotation
\[
Rot_s := \begin{pmatrix}\cos s  & -\sin s\\ \sin s& \cos s\end{pmatrix},
\]
and $\bar \bS_t$ is the unique solution of the $2$-D SDE
\[
d\bar \bS_t = -\bar \bS_t dt + d \B_t,
\]
with the initial value $\bar \bS_0=\frac{\sqrt{\mu}}{\sigma}\bar \Y_0$.
Therefore, $\|\bar \Y_t\|$ can be approximated by $R_t := \|\bar \Y^{app}_t\| = \frac{\sigma}{\sqrt{\mu}}\|\bar \bS_{\mu t}\|$
which by  Ito calculus satisfies the SDE 
\begin{equation}\label{CIRmodified}
\begin{split}
dR_t &=  \Big[\frac{\sigma^2}{2R_t} - \mu R_t \Big] dt + \sigma d\tilde B_t.
\end{split}
\end{equation}

The second method is to consider $\bar{\Y}_t$ in  polar coordinates with
\[
d \bar{\Y}_t = A \bar{\Y}_t dt + \h_e dB_t,
\]
where $\h_e=Q^{-1}\begin{pmatrix}0\\ \sigma_0\end{pmatrix} $. Its norm $\bar{R}_t:= \|\bar{\Y_t}\|$ and its angle $\thetab_t = \frac{\bar{\Y}_t}{\bar{R}_t}$ satisfy
\begin{eqnarray*}
    d\bar{R}_t &=&\Big[\frac{\|\h_e\|^2-\langle \h_e,\thetab_t\rangle^2}{2 \bar{R}_t} - \mu \bar{R}_t \Big] dt + \langle \thetab_t,\h_e\rangle dB_t, \\
    d\thetab_t &=& \Big[(A+\mu I) \thetab_t - \frac{\|\h_e\|^2 - \langle \h_e,\thetab_t \rangle^2}{2 \bar{R}_t^2} \thetab_t\Big] dt + \frac{1}{\bar{R}_t}\Big[ \h_e - \langle \h_e, \thetab_t \rangle \thetab_t\Big]dB_t.
\end{eqnarray*}
By the averaging technique from \cite[Theorem 1]{Bax11}, one can approximate $\thetab_t = \left(
\begin{array}{c}
 \sin \nu t \\
\cos \nu t \\
\end{array}
\right)$, hence
\begin{equation}\label{CIRmodified1}
\begin{split}
d\bar{R}_t
&=\Big[\frac{157.881  \sigma_0^2 - ( 1.27722 \sin \nu t  +  12.5 \cos \nu t )^2 \sigma_0^2}{2 \bar{R}_t} - \mu \bar{R}_t \Big] dt \\
&\qquad + ( 1.27722 \sin \nu t  +  12.5 \cos \nu t ) \sigma_0 dB_t.
\end{split}
\end{equation}

\begin{figure}[h!]
\begin{center}
    \includegraphics[width=8cm]{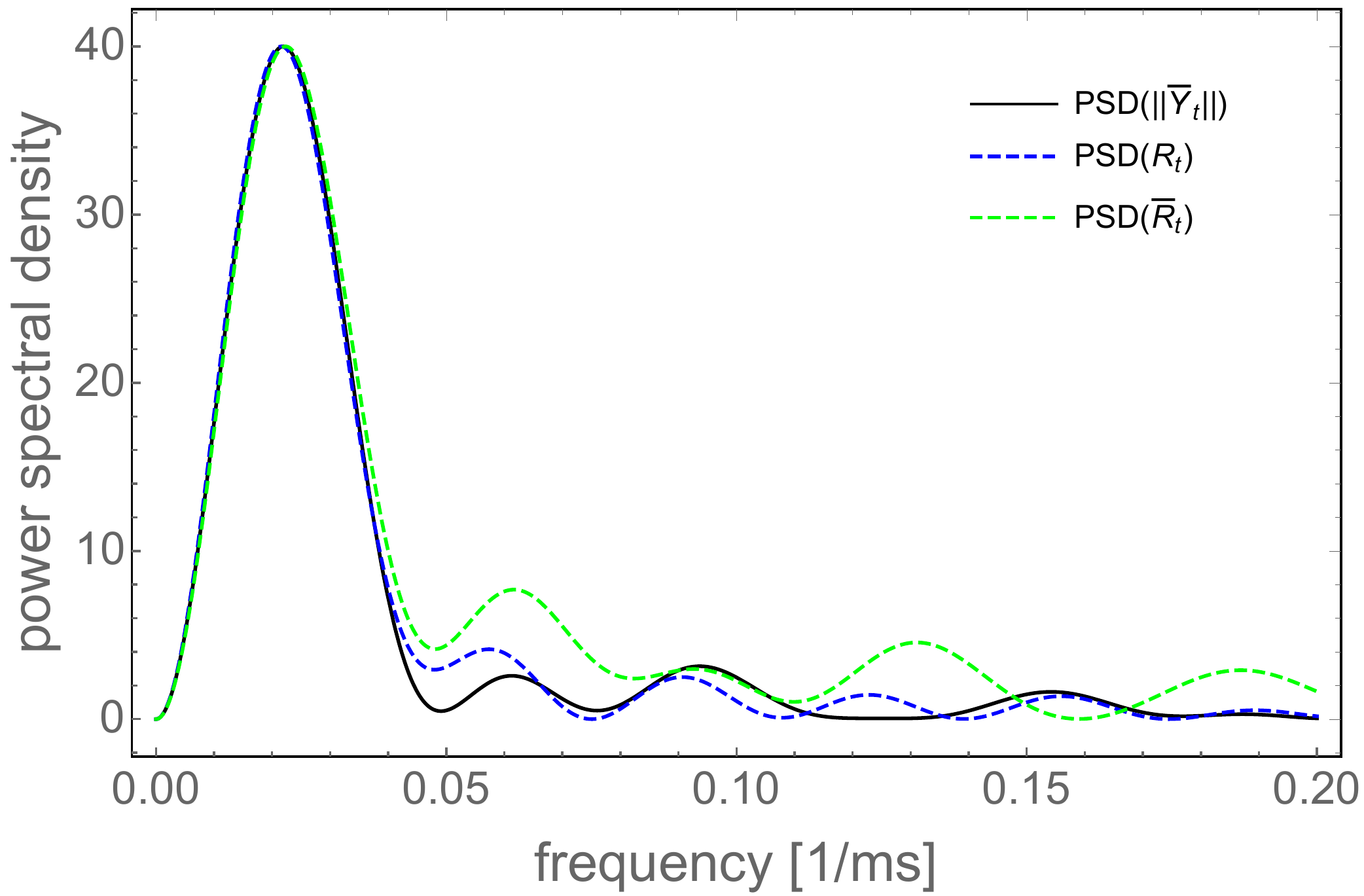}
\end{center}
\caption{\small The norm of the process $\{\bar{\Y}_t\}$ \eqref{additive1}, the process $\{R_t\}_t$~\eqref{CIRmodified}, and the process $\{\bar{R}_t\}_t$~\eqref{CIRmodified1} with the chosen parameters of system, $\sigma_0=0.01$ and the same starting point $\|Q^{-1} \bar \X_0\|$ are compared by using the power spectral density. Their spectrum densities are well approximated.}
\label{Fig:51}
\end{figure}

Thus, by using the averaging technique, we proved that both
Eqs.~\eqref{CIRmodified} and \eqref{CIRmodified1} are good
approximations of the radial process $\{\|\bar \Y_t\|\}_t =
\{\|Q^{-1}\bar \X_t\|\}_t$. This can also be tested by using the
power spectral density estimation (see Fig.~\ref{Fig:51}).

\subsection*{Firing mechanism}
A spike in Eq.~\eqref{eq:1} occurs when there is a transition of a
random trajectory from the vicinity  of the stable fixed point
$\X_e=(v_e,w_e)$ located on the left stable part of $\mathcal C_0$ to
its right stable part and back to the vicinity of $\X_e$.
This spike happens almost surely when a random trajectory with
the starting point $\X_0$ in the vicinity of $\X_e$ crosses the threshold
line $v=0$. From the phase space of
Eq.~\eqref{eq:1} (see Fig.~\ref{Fig:2}), the probability of a spike increases as the starting point $\X_0$ moves farther away from $\X_e$.

In order to construct the firing mechanism of Eq.~\eqref{CIRmodified}
matching that of Eq.~\eqref{eq:1}, we will calculate the conditional probability
that Eq.~\eqref{eq:1} fires given that the trajectory crosses the line
$L=\{(v_e,w): w \le w_e\}$. Denote by $L_i = (v_e,w_e-l_i)$ with $l_i = i\delta = i \frac{|w_e+0.453|}{20}$ for $i=0,1,\ldots,34$,
then the distance between the equilibrium and $L_i$ is $l_i$. The value $|w_e+0.453|$ can be considered as the distance between the fixed point $(v_e,w_e)$ and the separatrix (see also Fig.~\ref{Fig:1}) along $L$. For a given pair ($\sigma_0, l_i$), a short trajectory starting in $L_i$ was simulated from~\eqref{eq:1}, it was recorded whether a spike occurred (crossing the threshold $v=0$) in the first cycle of the stochastic path around $(v_e,w_e)$. This was repeated 1000 times and we counted the ratio of the number of spikes, denoted by $\hat p(l_i,\sigma_0)$, which is an estimate for the conditional probability of firing $p(l,\sigma_0)$. The estimation was, furthermore, repeated for $\sigma_0 = 0.001, 0.002, 0.003, 0.004, 0.005, 0.006, 0.007, 0.008, 0.009, 0.01, 0.015$.

\begin{figure}[h!]
\begin{center}
\includegraphics[width=7cm]{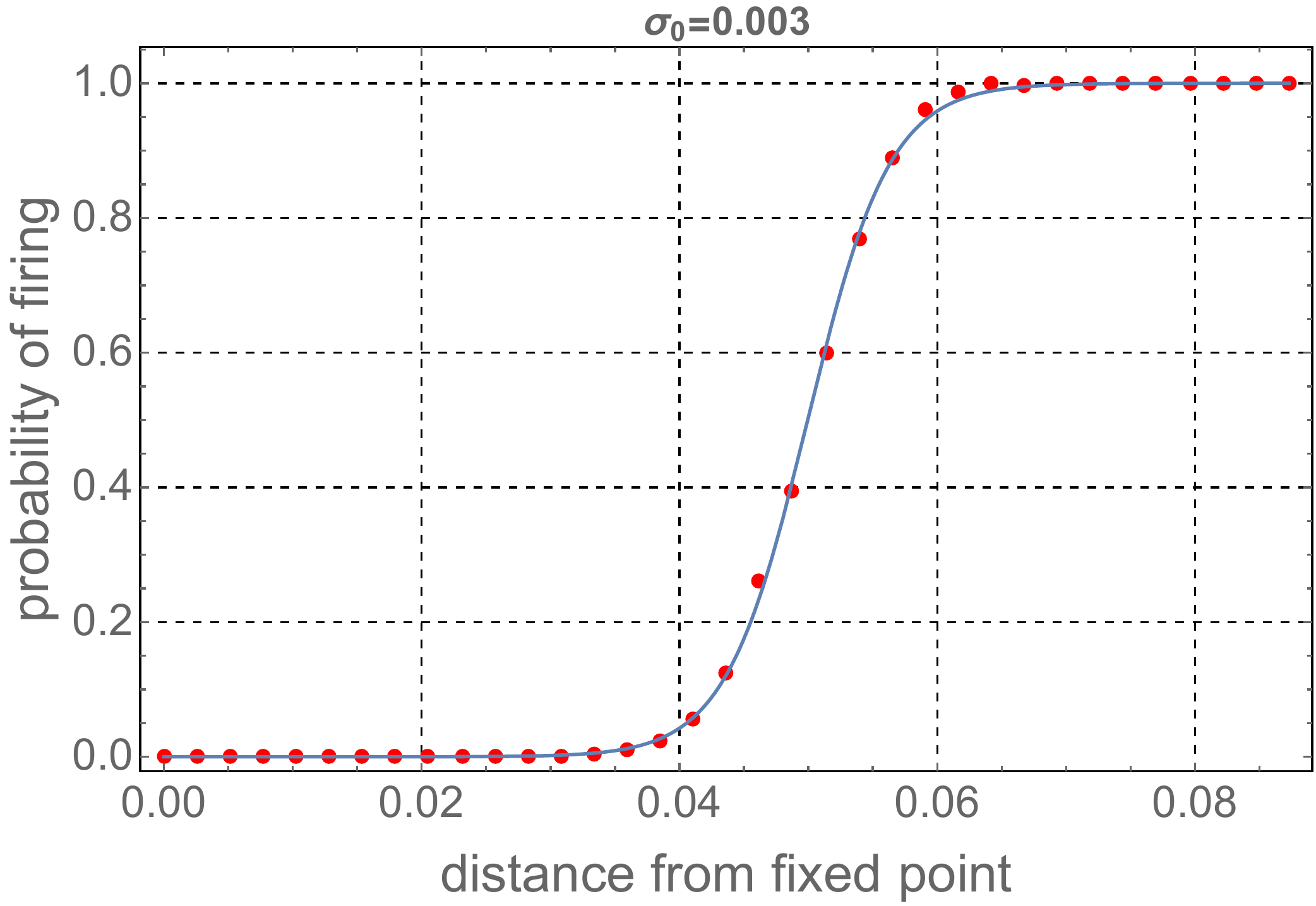}\qquad \qquad
\includegraphics[width=7cm]{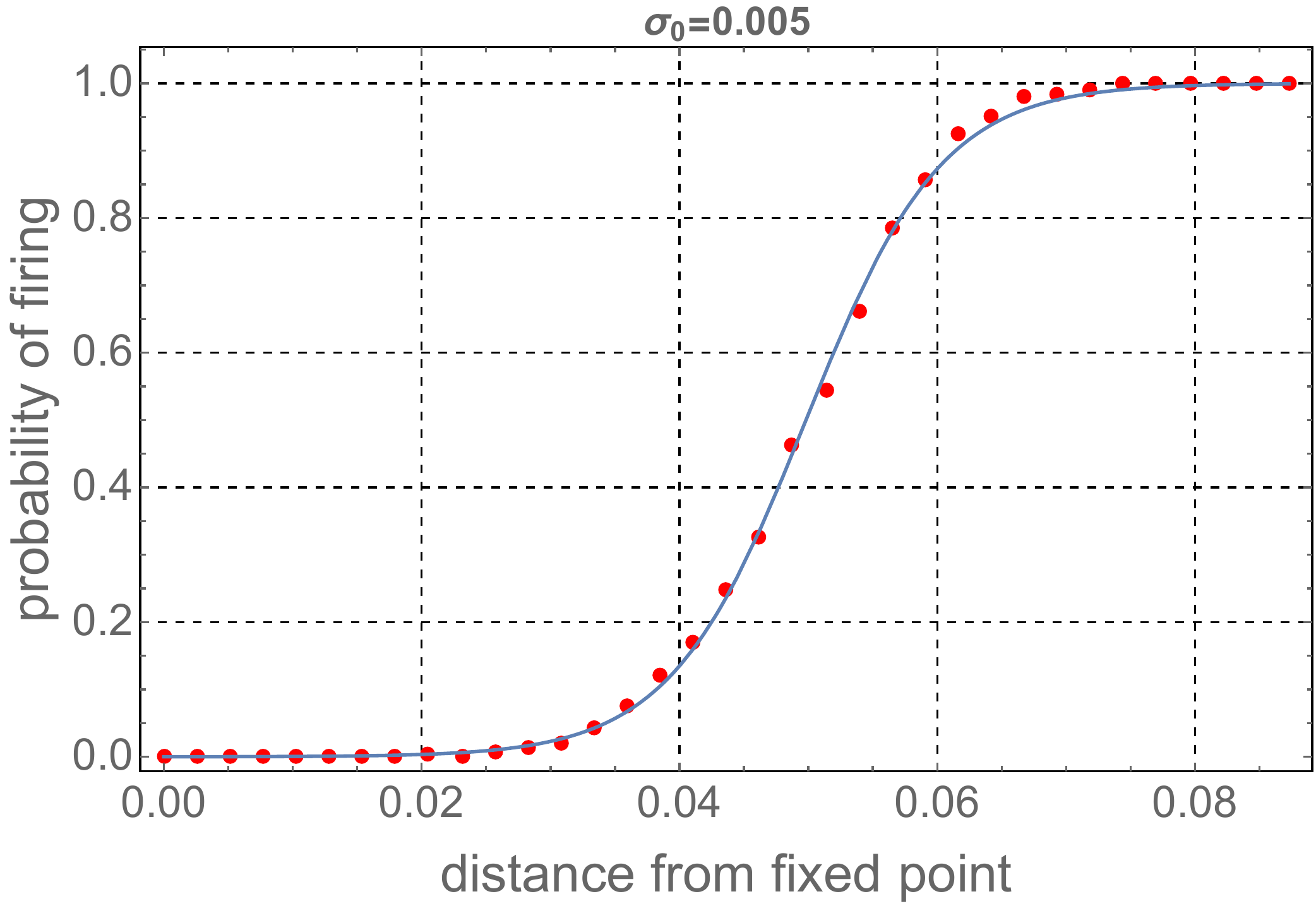}\\
\includegraphics[width=7cm]{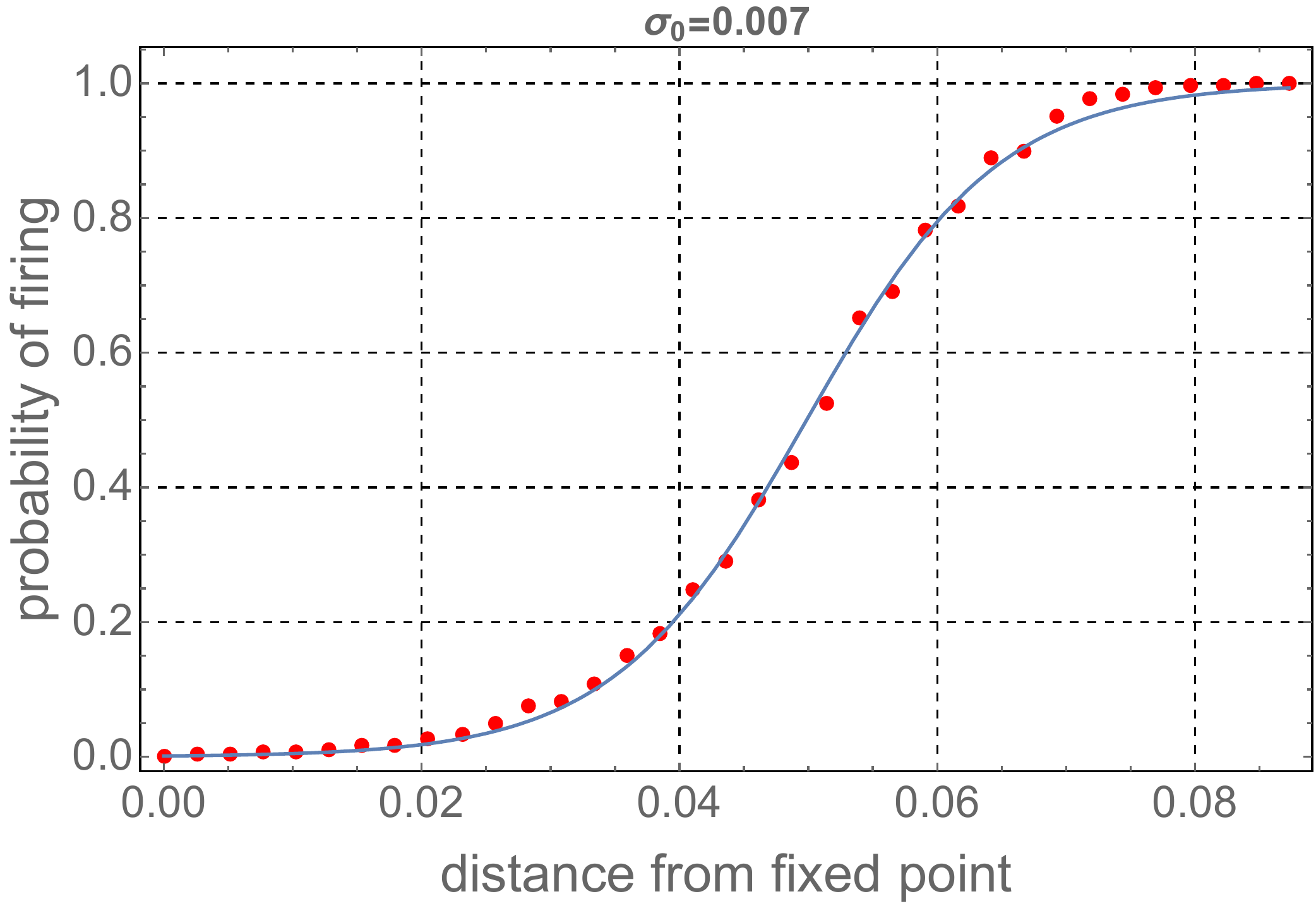}\qquad\qquad
\includegraphics[width=7cm]{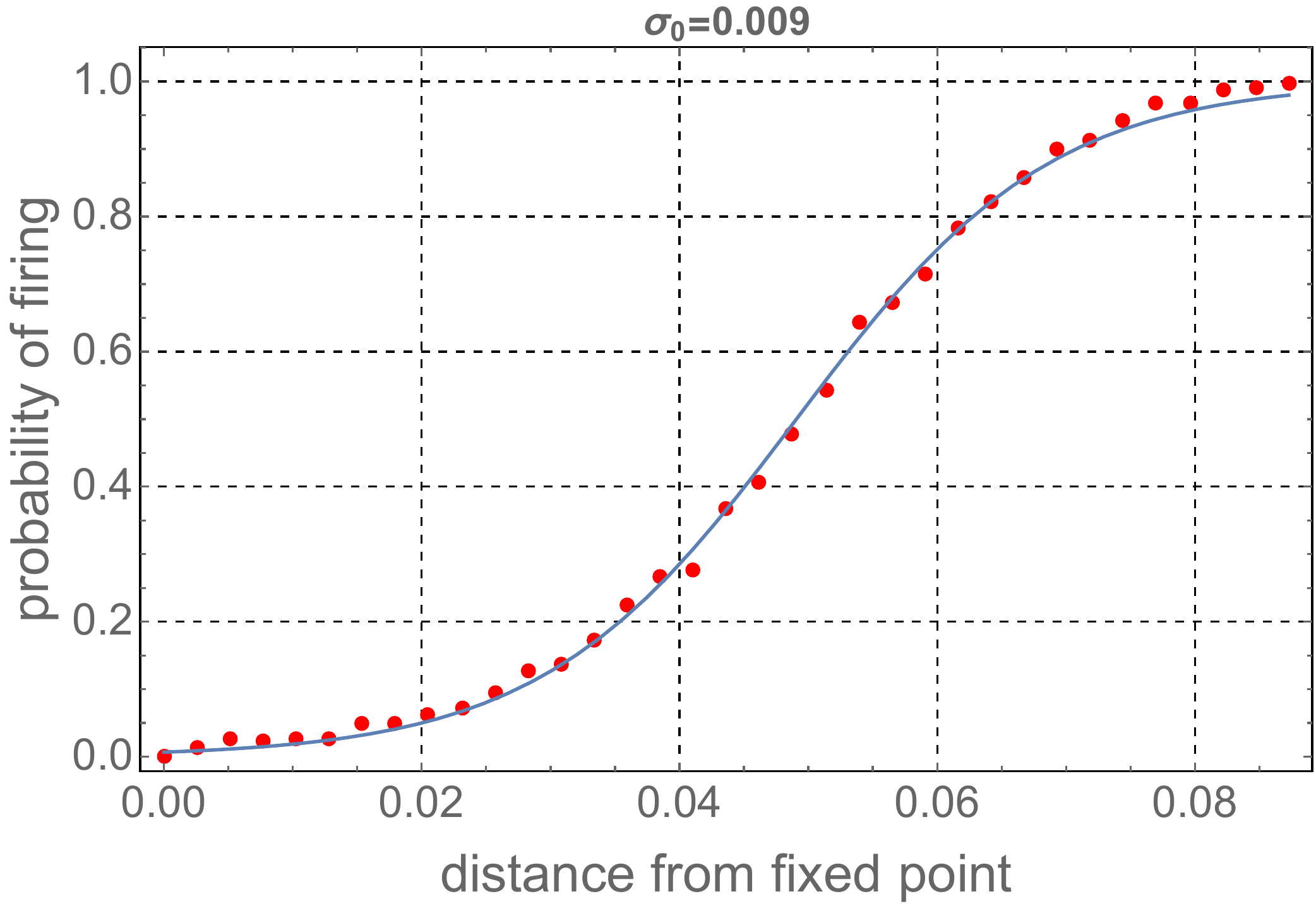}\\
\includegraphics[width=7cm]{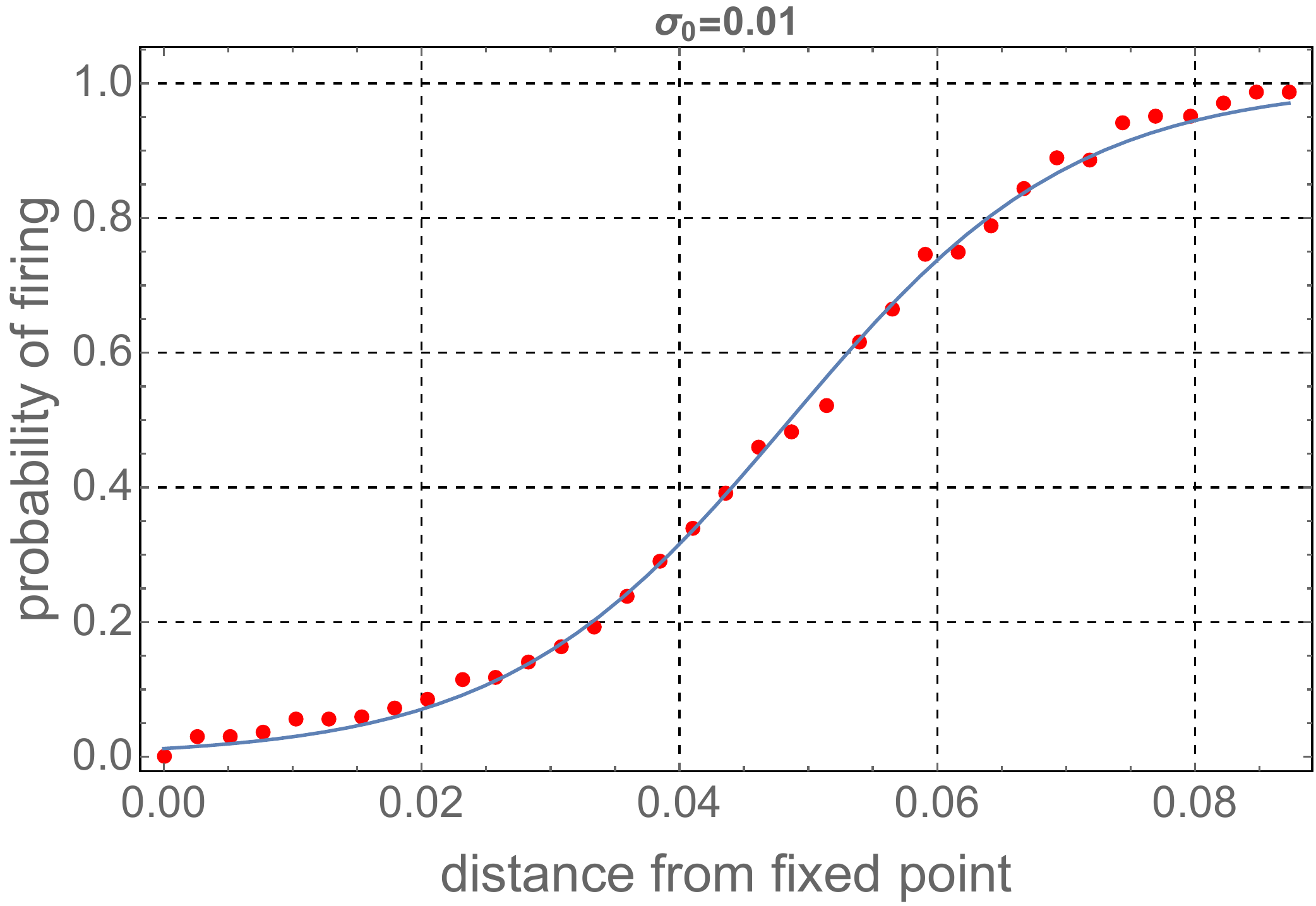} \qquad\qquad
\includegraphics[width=7cm]{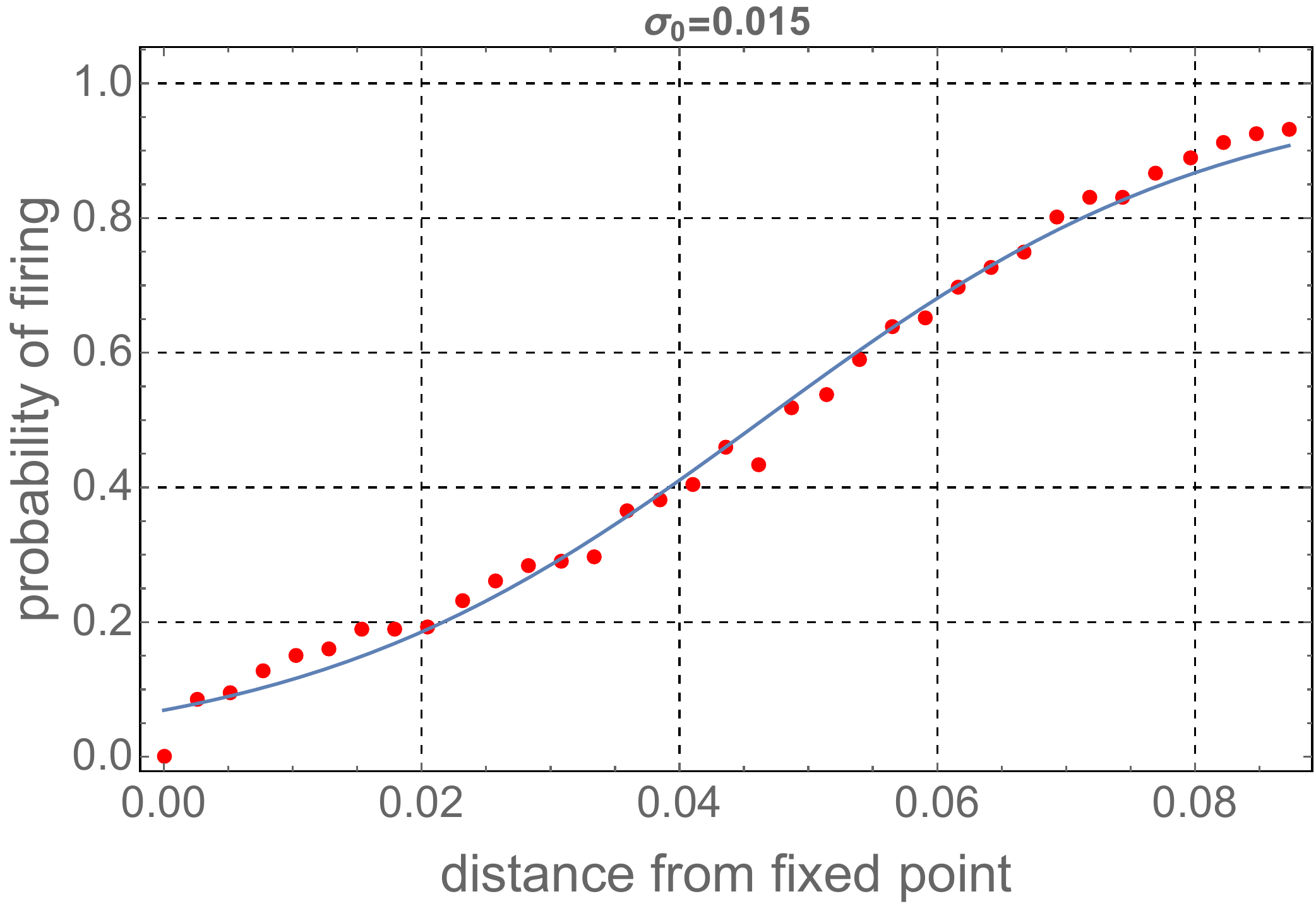}\\
\caption{Conditional probability of spiking when crossing the line $L=\{(v_e,w): w \le w_e\}$ for different values of the noise amplitude $\sigma_0$. The red dots are individual nonparametric estimates and the blue curve are the fitted curves given by~\eqref{eq:conditionalprob}.}\label{Fig:6}
\end{center}
\end{figure}

From the numerical simulation, for each $\sigma_0$, the estimate of the conditional probability is close to zero when we start in the
immediate neighborhood of the stable fixed point and close to one when we start at
the $L_{34}$, i.e., sufficiently far from the fixed point. Theses estimates appear to depend in a sigmoidal way on the distance from the stable fixed point. Therefore we assumed  the conditional probability of firing to be of the form
\begin{equation}\label{eq:conditionalprob}
p(l)= \frac{1}{1+e^{\frac{a-l}{b}}}.
\end{equation}
The parameters $a$ and $b$ then are estimated by using a non-linear regression from the above simulation data and are plotted in Fig.~\ref{Fig:6}  for some different values of the noise amplitude $\sigma_0 = 0.003, 0.005, 0.007, 0.009, 0.01$, and $0.015$. We see that the family of estimates, $\hat p$, fits the fitted curve quite well for each value of $\sigma_0$. Regression estimates are reported in Table~\ref{t:1}. Note that $p(a)=1/2$, i.e., $a$ is the distance along $L$ from $w_e$ at which the conditional probability of firing equals one half. For all values of $\sigma_0$, the estimate of $a$ is close to the distance along $L$ between $w_e$ and the separatrix, which equals $0.05$. In other words, the probability of firing, if the path starts at the intersection of $L$ with the sepametrix, is about $1/2$. The estimate of $b$ increases with respect to $\sigma_0$, and the conditional probability approaches a step function as the amplitude of the noise goes to zero. A step function would correspond to the firing being represented by a first passage time of a fixed threshold.

{
\renewcommand{\arraystretch}{1.1}
\begin{table}[h!]
\begin{center}
\small\addtolength{\tabcolsep}{-3.5pt}
\caption{Estimates of regression parameters for the conditional probability of firing in the original space and in the
transformed coordinates based on the additive noise $\sigma_0$}
\begin{tabular}{l l l l l l l l l l l l}
\hline
$\sigma_0$ & 0.001 & 0.002 & 0.003 & 0.004 & 0.005 & 0.006 & 0.007 & 0.008 & 0.009 &  0.01 & 0.015\\
\hline
\hline
$a$ 	&	0.050161 	&	0.050268	&	0.049946	 	&0.049760	&	0.049816 		&0.050001	&0.049862 	&0.049411 	&0.049078	&0.048559 	&0.046142\\
$b$ 	& 	0.001028 	& 	0.002099 	& 	0.003192 		&0.004310 	&	0.005281 		&0.006459 	&0.007478 	&0.008844 	&0.009877	&0.011068  	&0.017722\\
\hline
$a^*$ &	0.630282 	& 	0.631624 	&	0.627576  	&0.625240 	&	0.625935 		&0.628262 	&0.626516 	&0.620859 	&0.616673  	&0.610148	&0.579777\\
$b^*$ & 	0.012918 	&	0.026372 	&	0.040106   	&0.054158 	&	0.066352 		&0.081158 	&0.093960 	&0.111127 	&0.124107  	&0.139075	&0.222676\\
\hline
\end{tabular}\label{t:1}
\end{center}
\end{table}
}

To simplify calculations we will work on the transformed coordinates $\bar \Y_t$. Then the distance $l$ between $ (0,l)$ and $(0,0)$ in $\bar \X_t$ transforms to the distance
\[
r = \Bigg|  Q^{-1} \begin{pmatrix}0  \\l\end{pmatrix} \Bigg| = \sqrt{-\frac{m_{12}}{m_{21} \nu ^2}} l.
\]
and the conditional probability of firing Eq.~\eqref{eq:conditionalprob} transforms to
\begin{equation}\label{eq:conditionalprobTransformed}
p(r)= \frac{1}{1+e^{\frac{a^*-r}{b^*}}}
\end{equation}
where $a^* = \sqrt{-\frac{m_{12}}{m_{21} \nu ^2}}  a $ and $b^* = \sqrt{-\frac{m_{12}}{m_{21} \nu ^2}}  b$.

\subsection*{ISI distributions} 
\begin{figure}[h!]
\begin{center}
    \includegraphics[width=7cm]{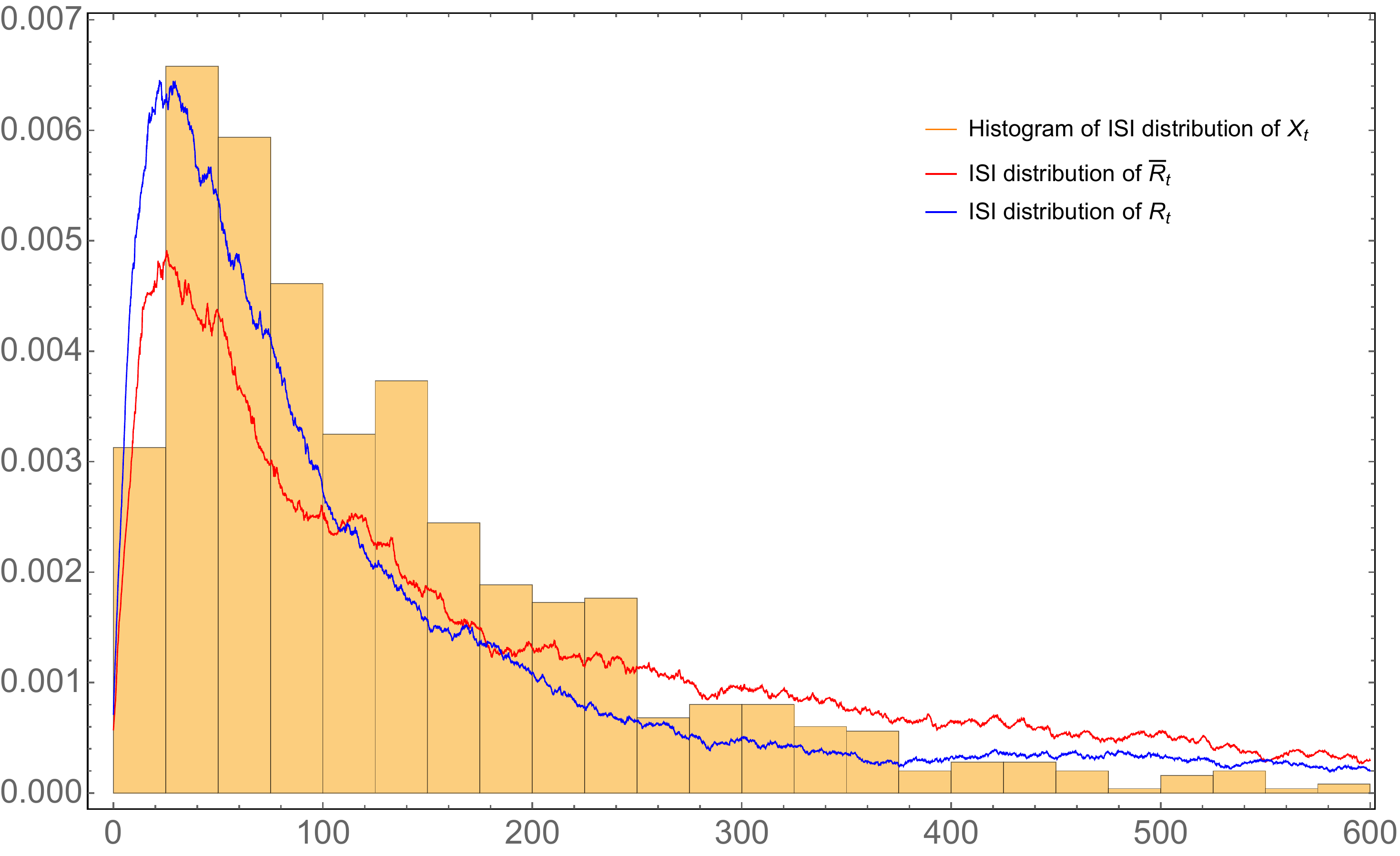}
\end{center}
\caption{\small The estimated ISI distributions from our approximate LIF models~\eqref{CIRmodified} and~\eqref{CIRmodified1} with the firing mechanism compare well with the estimated ISI histogram of FHN~\eqref{eq:1} reset to $0$ after firings. $\sigma_0 = 0.01, M =1000, n=10.$}
\label{Fig:7}
\end{figure}

The comparison of the original stochastic FHN model \eqref{eq:1} and the two LIF models \eqref{CIRmodified} and \eqref{CIRmodified1} can be performed by studying the ISI statistics. Namely, one first simulates the trajectories of the system~\eqref{eq:1} with starting points $\X_0$ close to the fixed point $\X_e$ until the first spiking time, and thereafter resets to the starting points.
Due to Theorem~\ref{approximation}, we can simplify the simulation by choosing the starting point at exactly $\X_e$. This was done 1000 times, and the time of the first firing was recorded. A histogram for this data is shown in Fig.~\ref{Fig:7}. The ISI-distribution of Eq.~\eqref{CIRmodified}
is computed as follows (the ISI-distribution of
Eq.~\eqref{CIRmodified1} is computed similarly). Let $\tau_1$ be
the first firing time. We computed the density of the distribution
of $\tau_1$ in terms of the conditional hazard rate \cite{Dil12},
\[
\alpha(r,t) = \lim_{\Delta t \to 0} \frac{1}{\Delta t} P(t\le \tau_1 < t+\Delta t | \tau_1\ge t, R_{t} =r).
\]
This function is the density of the conditional probability, given the position on $L$ is $r$ at time $t$, of a spike occurring in the next small time interval, given that it has not yet occurred.\\ 
Notice that the estimated conditional probability of firing~\eqref{eq:conditionalprobTransformed} is calculated in one cycle of the process, which on average takes $2\pi/\nu$ time units. Therefore, we estimate
the hazard rate as
\begin{equation}\label{eq:hazardrate}
\alpha(r,t) = \alpha(r) = \frac{\nu}{2 \pi} \frac{1}{1+e^{\frac{a^*-r}{b^*}}} .
\end{equation}

On the other hand, from standard results from survival analysis, see e.g. \cite{aalen10} we know that the density of the firing time can be calculated as
\begin{equation}\label{eq:density}
g(t) = \frac{d}{dt} P(\tau_1 \le t) = E \Bigg( \alpha(R_t) e^{-\int_0^t \alpha(R_s)ds}\Bigg).
\end{equation}

Due to the law of large numbers, for fixed $t$, we can numerically determine the density~\eqref{eq:density} up to any desired precision by choosing $n$ and $M$ large enough through the expression
\[
g(t) \approx \frac{1}{M} \sum\limits_{m=1}^M \alpha(R^{(m)}_t) e^{-\frac{t}{n} \sum\limits_{i=1}^n \frac{\alpha\big(R^{(m)}_{it/n}\big)+\alpha\big(R^{(m)}_{(i-1)t/n}\big)}{2}}. 
\]
Here ($R^{(m)}_0, \ldots, R^{(m)}_{it/n}, \ldots, R^{(m)}_{t})$ are $M$ realizations of $R_{it/n}, i = 0,1,\ldots, n$, and the integral has been approximated by the trapezoidal rule. The results are illustrated in Fig.~\ref{Fig:7} for $\sigma_0=0.01$, using $M = 1000, n=10$. The estimated ISI distributions from our approximate LIF models~\eqref{CIRmodified} and~\eqref{CIRmodified1} with the firing mechanism compare well with the estimated ISI histogram of FHN~\eqref{eq:1} reset to $0$ after firings.

\section{Appendix}\label{section5}
\begin{proof}[{\bf Proof of Theorem \ref{stadis}}]
	We are going to prove that there exists a random pullback attractor for the general equation \eqref{SDE}. Consider two cases:\\
	
	\begin{itemize}
		\item Additive noise: In this case, the proof follows  similar steps as in  \cite{GA09}. We define $\Y_t = \X_t - \eta_t$ where $\eta_t$ is the unique stationary solution of
		\begin{eqnarray*}
			d\eta_t = -\eta_t dt + (0,\sigma_0)^T dB_t.
		\end{eqnarray*}
		System \eqref{SDE} is then tranformed to
		\begin{equation}\label{RDE2}
			\dot{\Y}_t = F(\Y_t + \eta_t) + \eta_t.
		\end{equation}
		Observe that
		\begin{eqnarray*}
			\frac{d}{dt} \|\Y_t\|^2
			&=& 2 \langle \Y_t, F(\Y_t + \eta_t) - F(\eta_t)\rangle + 2 \langle \Y_t, F(\eta_t) + \eta_t \rangle \\
			&\leq& 2(a-b\|\Y_t\|^2) + b\|\Y_t\|^2 + \frac{1}{b} \|F(\eta_t) + \eta_t\|^2 \\
			&=& 2a + \frac{1}{b} \|F(\eta_t) + \eta_t\|^2 - b \|\Y_t\|^2.
		\end{eqnarray*}
		Hence by the comparison principle, $\|\Y_t\| \leq R_t$ whenever $\|\Y_0\|^2 \leq R_0$ where $R_t$ is the solution of
		\begin{equation}\label{addReq}
			\dot{R}_t = 2a  + \frac{1}{b} \|F(\eta_t) + \eta_t\|^2 - b R_t,
		\end{equation}
		which can be computed explicitly as
		\[
		R_t(\omega,R_0) = e^{-bt}R_0 + \int_0^t e^{-b(t-s)} \Big[2a  + \frac{1}{b} \|F(\eta_s) + \eta_s\|^2 \Big] ds.
		\]
		It is then easy to check that the vector field in \eqref{RDE2} satisfies the local Lipschitz property and the solution is bounded and thus of linear growth on any fixed $[0,T]$, see e.g. \cite{hoppe96}. Hence there exists a unique solution of \eqref{RDE2} with initial condition, which also proves the existence and uniqueness of the solution of \eqref{SDE}. The cocycle property \eqref{cocycle}  follows automatically from \cite[Chapter 2]{arnold}.
		
		A direct computation shows that there exists a random radius
		\[
		R^*(\omega) = \int_{-\infty}^{0} [2a + \frac{1}{b} \|F(\eta_s) + \eta_s\|^2] e^{bs} ds,
		\]
		which is the stationary solution of \eqref{addReq}, such that $\X_t(\omega,\X_0) \in B(\eta_t,R^*(\theta_t \omega))$ whenever $\X_0 \in B(\eta_0,R^*(\omega))$ by the comparison principle, and furthermore,
		\[
		\limsup \limits_{t\to \infty}\|\Y_t(\theta_{-t}\omega,\Y_0)\|^2 \leq \limsup \limits_{t \to \infty} R_t(\theta_{-t} \omega,R_0) = R^*(\omega).
		\]
		Hence the random ball $B(\eta,R^*)$ is a forward invariant pullback absorbing set of the random dynamical system generated by $\varphi(t,\omega)\X_0$ \eqref{SDE}. By the classical theorem \cite{crauel}, there exists the global random pullback attractor for \eqref{SDE}. \\
		
		\item Multiplicative noise:  In this case, we introduce the transformation
		\begin{eqnarray}\label{conjugacy}
			\Y_t = (v_t,\bar{\omega}_t)^T:= \begin{pmatrix}
				1&0\\0&e^{-\sigma_0 z_t}
			\end{pmatrix}\X_t = T(z_t)\X_t
		\end{eqnarray}
		where $z_t$ is the unique stationary solution of the Ornstein-Uhlenbeck equation
		\begin{equation}\label{Langevineq}
			d z_t = - z_t dt + dB_t.
		\end{equation}
		This transforms system \eqref{SDE}  into a random differential equation.
		\begin{eqnarray}\label{RDE}
			\dot{v_t} &=& v_t -\frac{v_t^3}{3} - e^{\sigma_0 z_t} \bar{\omega}_t + I \\
			\dot{\bar{\omega}}_t &=& e^{-\sigma_0 z_t} \varepsilon v_t + (\sigma_0 z_t - \varepsilon \beta) \bar{\omega}_t + \varepsilon \alpha e^{-\sigma_0 z_t}.\notag
		\end{eqnarray}
		or equivalently,
		\[
		\dot{\Y}_t = G(z_t,\Y_t)
		\]
		where $G$ satisfies $G(z_t,0) = (I,\varepsilon \alpha e^{-\sigma_0 z_t})^T$ and
		\begin{eqnarray*}
			&&\langle \Y_1 - \Y_2,G(z_t,\Y_1)-G(z_t,\Y_2) \rangle \\
			&=& (v_1-v_2)^2 \Big[1- \frac{1}{3} (v_1^2 + v_1 v_2 + v_2^2)\Big]+(\epsilon e^{-\sigma_0 z_t} - e^{\sigma_0 z_t})(v_1-v_2)(\bar{w}_1-\bar{w}_2)\\
			&& +(\sigma_0 z_t -\epsilon \beta) (\bar{w}_1-\bar{w}_2)^2\notag\\
			&\leq& (v_1-v_2)^2 - \frac{1}{12}(v_1-v_2)^4 + \frac{1}{2\epsilon \beta} (\epsilon e^{-\sigma_0 z_t} - e^{\sigma_0 z_t})^2 (v_1-v_2)^2 \\
			&& + (\sigma_0 z_t -\frac{\epsilon \beta}{2}) (\bar{w}_1-\bar{w}_2)^2\notag\\
			&\leq& - \frac{1}{12}(v_1-v_2)^4 + \Big[1+ \frac{1}{2\epsilon \beta} (\epsilon e^{-\sigma_0 z_t} - e^{\sigma_0 z_t})^2 - \sigma_0 z_t +\frac{\epsilon \beta}{2}\Big](v_1-v_2)^2 \\
			&& + (\sigma_0 z_t -\frac{\epsilon \beta}{2}) \|\Y_1-\Y_2\|^2\notag\\
			&\leq& -\frac{1}{12} \Big((v_1-v_2)^2  + 6 \Big[1+ \frac{1}{2\epsilon \beta} (\epsilon e^{-\sigma_0 z_t} - e^{\sigma_0 z_t})^2 - \sigma_0 z_t +\frac{\epsilon \beta}{2}\Big] \Big)^2 \\
			&&+ 3 \Big[1+ \frac{1}{2\epsilon \beta} (\epsilon e^{-\sigma_0 z_t} - e^{\sigma_0 z_t})^2 - \sigma_0 z_t +\frac{\epsilon \beta}{2}\Big]^2 + (\sigma_0 z_t -\frac{\epsilon \beta}{2}) \|\Y_1-\Y_2\|^2\notag\\
			&\leq& 3 \Big[1+ \frac{1}{2\epsilon \beta} (\epsilon e^{-\sigma_0 z_t} - e^{\sigma_0 z_t})^2 - \sigma_0 z_t +\frac{\epsilon \beta}{2}\Big]^2 + (\sigma_0 z_t -\frac{\epsilon \beta}{2}) \|\Y_1-\Y_2\|^2.
		\end{eqnarray*}
		Thus,
		\begin{eqnarray*}
			\frac{d}{dt} \|\Y_t\|^2 &=& 2 \langle \Y_t-0, G(z_t, \Y_t) - G(z_t,0)\rangle + 2 \langle \Y_t, G(z_t,0) \rangle\\
			&\leq&3 \Big[1+ \frac{1}{2\epsilon \beta} (\epsilon e^{-\sigma_0 z_t} - e^{\sigma_0 z_t})^2 - \sigma_0 z_t +\frac{\epsilon \beta}{2} \Big]^2 +(\sigma_0 z_t -\frac{\epsilon \beta}{2}) \|\Y_t\|^2 \\
			&& + 2 \langle \Y_t, G(z_t,0) \rangle\\
			&\leq& 3 \Big[1+ \frac{1}{2\epsilon \beta} (\epsilon e^{-\sigma_0 z_t} - e^{\sigma_0 z_t})^2 - \sigma_0 z_t +\frac{\epsilon \beta}{2} \Big]^2 + \frac{4}{\varepsilon \beta} \|G(z_t,0)\|^2\\
			&& +(\sigma_0 z_t -\frac{\epsilon \beta}{4}) \|\Y_t\|^2 \\
			&\leq&3 \Big[1+ \frac{1}{2\epsilon \beta} (\epsilon e^{-\sigma_0 z_t} - e^{\sigma_0 z_t})^2 - \sigma_0 z_t +\frac{\epsilon \beta}{2} \Big]^2 + \frac{4}{\varepsilon \beta} \Big[I^2 + \varepsilon^2 \alpha^2 e^{-2 \sigma_0 z_t}\Big]\\
			&&+(\sigma_0 z_t -\frac{\epsilon \beta}{4}) \|\Y_t\|^2 \\
			&\leq&  p(z_t) + q(z_t) \|\Y_t\|^2.
		\end{eqnarray*}
		Hence by the comparison principle, $\|\Y_t\|^2 \leq R_t$ whenever $\|\Y_0\|^2 \leq R_0$ where $R_t$ is the solution of
		\begin{equation}\label{compared}
			\dot{R}_t = p(z_t) + q(z_t) R_t,
		\end{equation}
		which can be computed explicitly as
		\[
		R_t(\omega,R_0) = e^{\int_0^t q(z_u(\omega)) du} R_0 + \int_{0}^{t} p(z_s(\omega)) e^{\int_s^t q(z_u(\omega)) du} ds.
		\]
		Using similar arguments as in the additive noise case, there exists a unique solution of \eqref{RDE} and \eqref{SDE}. Also, the solution generates a random dynamical system.
		
		On the other hand, observe that by the Birkhorff ergodic theorem, there exists almost surely
		\[
		\lim \limits_{t \to -\infty}\frac{1}{t} \int_t^0 q(z_u) du = \lim \limits_{t\to -\infty}\frac{1}{t} \int_t^0 q(z(\theta_u \omega)) = E \Big[\sigma_0 z(\cdot) - \frac{\varepsilon \beta}{4}\Big] = -\frac{\varepsilon \beta}{4} <0,
		\]
		therefore there exists a unique stationary solution of \eqref{compared} which can be written in the form
		\[
		\bar{R}(\omega) = \int_{-\infty}^{0} p(z_s(\omega)) e^{\int_s^0 q(z_u(\omega)) du} ds.
		\]
		Moreover, $\|\Y_t(\omega,\Y_0)\|^2 \leq \bar{R}(\theta_t \omega)$ whenever $\|\Y_0\|^2 \leq \bar{R}(\omega)$ and
		\[
		\limsup \limits_{t\to \infty}\|\Y_t(\theta_{-t}\omega,\Y_0)\|^2 \leq \limsup \limits_{t \to \infty} R_t(\theta_{-t} \omega,R_0) = \bar{R}(\omega).
		\]
		Hence, the ball $B(0,R(\omega))$ is actually forward invariant under the random dynamical system generated by \eqref{RDE} and is also a pullback absorbing set. Again by applying \cite{crauel}, there exists a random attractor for \eqref{RDE}. Due to the fact that $z_t$ is the stationary solution of \eqref{Langevineq}, it is easy to see that the random linear transformation $T(z)$ given in \eqref{conjugacy} is {\it tempered} (see \cite[pp. 164, 386]{arnold}), i.e.
		\[
		0\leq \lim \limits_{t \to \infty} \frac{1}{t} \log \|T(z_t)\| = \lim \limits_{t \to \infty} \frac{1}{2t} \log (1+ e^{-2\sigma_0 z_t}) \leq \lim \limits_{t \to \infty} \frac{1}{2t} (1 + 2\sigma_0 |z_t|) = 0.
		\]
		Therefore, it follows from \cite{imkeller} that systems \eqref{SDE} and \eqref{RDE} are conjugate under the tempered transformation \eqref{conjugacy}, hence there exists also a random attractor for system \eqref{SDE}.
	\end{itemize}
	
\end{proof}

\begin{proof}[{\bf Proof of Theorem \ref{approximation}}]
	Observe that the matrix
	\[
	DF(\X_e) = \begin{pmatrix} m_{11}&m_{12}\\ m_{21}& m_{22} \end{pmatrix}
	\]
	has two conjugate complex eigenvalues with negative real part
	\[
	\lambda_{1,2}  = \frac{1}{2}(1-v_e^2 - \epsilon \beta) \pm \frac{i}{2} \sqrt{4\epsilon - (1-v_e^2+\epsilon \beta)^2}=-0.0730077\pm 0.31615 i = -\mu \pm \nu i.
	\]
	Hence by using the transformation $\X-\X_e=Q \Y$ and $\bar{X} = Q \bar{\Y}$ with
	\[
	Q=\begin{pmatrix} -\nu&m_{11}+\mu \\ 0& m_{21} \end{pmatrix},
	\]
	the equations \eqref{shifteq} and \eqref{eq:O1} are transformed into the normal forms
	\begin{eqnarray}\label{transeq}
	d\Y_t &=& \Big[Q^{-1}DF(\X_e)Q \Y_t + Q^{-1}\bar{F}(Q \Y_t) \Big] dt + Q^{-1}H(Q\Y_t+ \X_e) \circ dB_t \\
	&=& [A \Y_t + F_1(\Y_t)] dt + Q^{-1}H(Q\Y_t+ \X_e) \circ dB_t,\\ \notag
	\Y_0 &=& Q^{-1}(\X_0-\X_e),   \notag
	\end{eqnarray}
	and
	\begin{eqnarray}\label{translin}
	d{\bar\Y}_t &=&  A \bar{\Y}_t  dt + Q^{-1}H(Q\bar{\Y}_t+\X_e) \circ dB_t,\\
	\bar{\Y}_0 &=& Q^{-1}(\X_0 -\X_e). \notag
	\end{eqnarray}
	where
	\[
	A = Q^{-1} DF(\X_e) Q = \begin{pmatrix}-\mu & \nu\\ -\nu & -\mu\end{pmatrix};\qquad F_1(\Y) :=Q^{-1}\bar{F}(Q\Y),
	\]
	and
	\begin{equation}\label{nlest}
	\| F_1(\Y)\|\leq \gamma(r) \|Q^{-1}\|\|Q\Y\| \leq \|Q^{-1}\|\gamma(r)r,\qquad \forall \|\Y\|\leq \frac{r}{\|Q\|}.
	\end{equation}
	Define the difference $\Z_t := \Y_t - \bar{\Y}_t$, then $\Z_t$ satisfies
	\begin{eqnarray*}
		d\Z_t &=& [A \Z_t + F_1(\Y_t)]dt + B_1 \Z_t \circ dB_t\\
		&=& \Big[(A+\frac{1}{2}B_1^TB_1) \Z_t + F_1(\Y_t)\Big]dt + B_1 \Z_t dB_t,
	\end{eqnarray*}
	where
	\[
	B_1 := 0 \text{\ if\ } H(\X) = (0,\sigma_0)^T\quad \text{and}\quad B_1 := Q^{-1}BQ \text{\ if\ } H(\X) = B\X.
	\]
	We analyze these two cases separately. \\
	
	\begin{itemize}
		\item Additive noise: then the equation for $\Z_t$ becomes deterministic, hence
		\begin{eqnarray*}
			\frac{d}{dt}\|\Z_{t\wedge\tau}\|^2 &=&2 \Big\langle \Z_{t\wedge\tau}, A \Z_{t\wedge\tau} + F_1(\Y_{t\wedge\tau})\Big \rangle \\
			&\leq& - 2 \mu \|\Z_{t\wedge\tau}\|^2 + \mu \|\Z_{t\wedge\tau}\|^2 + \frac{1}{\mu} \| F_1(\Y_{t\wedge\tau})\|^2 \\
			&\leq& \frac{1}{\mu} \|Q^{-1}\|^2 \gamma(r)^2 r^2 - \mu  \|\Z_{t\wedge\tau}\|^2.
		\end{eqnarray*}
		Using the fact that $\Z_0 = 0$, it follows that
		\[
		\|\Z_{t\wedge\tau}\|^2 \leq \frac{1}{\mu^2} \|Q^{-1}\|^2 \gamma(r)^2 r^2  + e^{-\mu (t\wedge \tau)} \Big(\|\Z_0\|^2 - \frac{1}{\mu^2} \|Q^{-1}\|^2 \gamma(r)^2 r^2 \Big).
		\]
		Therefore,
		\begin{equation*}
			\sup_{t \leq \tau} \|\Z_t\| \leq \frac{1}{\mu} \|Q^{-1}\| \gamma(r) r
		\end{equation*}
		which proves \eqref{estthm1} with $C = \frac{1}{\mu} \|Q\| \|Q^{-1}\|$.
		\item Multiplicative noise: By Ito's formula for the stopping time,
		\begin{eqnarray*}
			d\|\Z_{t\wedge\tau}\|^2 &=&2 \Big\langle \Z_{t\wedge\tau}, (A+\frac{1}{2}B_1^TB_1) \Z_{t\wedge\tau} + F_1(\Y_{t\wedge\tau})\Big \rangle d(t\wedge\tau) + \|B_1 \Z_{t\wedge\tau}\|^2 d({t\wedge\tau}) \\
			&&+ 2 \langle \Z_{t\wedge\tau},B_1 \Z_{t\wedge\tau}\rangle d B_{t\wedge\tau},
		\end{eqnarray*}
		hence taking the expectation on both sides and using \eqref{nlest} we have
		\begin{eqnarray*}
			\frac{d}{dt}E\|\Z_{t\wedge\tau}\|^2
			&\leq& 2\Big(-\mu + \|B_1^TB_1\|\Big)E\|\Z_{t\wedge\tau}\|^2 + 2 \|Q^{-1}\|\gamma(r)r E \|\Z_{t\wedge\tau}\|\\
			&\leq& (-\mu + 2 \|B_1^TB_1\|)E\|\Z_{t\wedge\tau}\|^2 \\
			&&+ \Big[-\mu \big(E\|\Z_{t\wedge\tau}\|\big)^2 +2 \|Q^{-1}\|\gamma(r)r E \|\Z_{t\wedge\tau}\|  \Big] \\
			&\leq& (-\mu + 2 \|B_1^TB_1\|)E\|\Z_{t\wedge\tau}\|^2 + \frac{1}{\mu} \|Q^{-1}\|^2\gamma(r)^2r^2,
		\end{eqnarray*}
		where the last inequality follows from the Cauchy inequality. Since
		\begin{equation}\label{sigma0}
		\lambda = \mu - 2 \|B^T_1 B_1\| >0,
		\end{equation}
		by noting that $\Z_0 = 0$, we get
		\begin{eqnarray*}
			E\|\Z_{t\wedge\tau}\|^2 &\leq& E\|\Z_{0\wedge\tau}\|^2 e^{-\lambda (t\wedge\tau)}+ \frac{1}{\mu} \|Q^{-1}\|^2\gamma(r)^2r^2 \ \frac{1}{\lambda} \Big[1 - e^{-\lambda(t\wedge\tau)}\Big]\\
			&\leq& \frac{1}{\mu} \frac{1}{\lambda}\|Q^{-1}\|^2\gamma(r)^2r^2,
		\end{eqnarray*}
		which proves \eqref{estthm} by choosing $C :=  \frac{1}{\mu} \frac{1}{\lambda}\|Q^{-1}\|^2 \|Q\|^2$.
	\end{itemize}
		
\end{proof}

\section*{Acknowledgments}
We thank the anonymous reviewers for their careful reading and useful remarks which helped to improve the quality of the manuscript. 

\bibliographystyle{apalike}

\end{document}